\documentclass[10 pt]{article}
\usepackage[pagebackref]{hyperref}
\hypersetup{
colorlinks = true,
linkcolor = blue,
citecolor = blue}
\usepackage[utf8]{inputenc}
\usepackage[T1]{fontenc}

\usepackage[]{amsfonts}
\usepackage{amssymb}
\usepackage{amsmath,amsthm}
\def\couleur(#1 #2 #3)
	{
\special{ps::[end]}}

\def\bx#1{\setbox1=\hbox{\kern3pt{#1}\kern3pt}			
 \dimen1=\ht1 \advance\dimen1 by 3pt \dimen2=\dp1 \advance\dimen2 by 3pt
 \setbox1=\hbox{\vrule height\dimen1 depth\dimen2\box1\vrule}%
 \setbox1=\vbox{\hrule\box1\hrule}%
 \advance\dimen1 by .4pt \ht1=\dimen1
 \advance\dimen2 by .4pt \dp1=\dimen2 \box1\relax}

\def\wbb#1{\kern#1em}
\def\vci{\vrule  width.02em height1.47ex depth-.0ex}		
\def\11{{\rm\wbb{.2}\vci\wbb{-.37}1}}

\def\underset#1#2{\mathrel{\mathop{\kern0pt #2}\limits_{#1}}}

\def\overset#1#2{\mathrel{\mathop{\kern0pt #2}\limits^{#1}}}

\newtheorem{thm}{Theorem}[section]
\newtheorem{lem}[thm]{Lemma}
\newtheorem{prop}[thm]{Proposition}
\newtheorem{cor}[thm]{Corollary}
\newtheorem{defin}[thm]{Definition}
\newtheorem{rem}[thm]{Remark}

\begin{document}

\title{Sobolev embeddings with weights in complete riemannian manifolds.}

\author{Eric Amar\protect\footnote{Univ. Bordeaux, CNRS, Bordeaux INP, IMB, UMR 5251,  F-33400, Talence, France}}

\date{ }
\maketitle
 \renewcommand{\abstractname}{Abstract}

\begin{abstract}
We prove Sobolev embedding Theorems with weights for vector bundles
 in a complete riemannian manifold. We also get general Gaffney's
 inequality with weights. As a consequence, under a "weak bounded
 geometry" hypothesis, we improve classical Sobolev embedding
 Theorems for vector bundles in a complete riemannian manifold.
 We also improve known results on Gaffney's inequality in a complete
 riemannian manifold.\par 

\end{abstract}
\ \par 

\tableofcontents
\ \par 
\ \par 

\section{Introduction.}
\quad Let $(M,g)$ be a complete riemann manifold. The Sobolev inequalities
 in $M$ for functions play a major role in the study
 of differential
 operators and nonlinear functional analysis. They are valid
 in ${\mathbb{R}}^{n}$ or if $M$ is compact. More generally we have:\ \par 

\begin{thm}
~\label{SC18}Let $(M,g)$ be a complete riemannian manifold of
 dimension $n$ with Ricci curvature bounded from below. The Sobolev
 embeddings for functions are valid for $M$ if and only if there
 is a uniform lower bound for the volume of balls which is independent
 of their center, namely if and only if $\inf _{x\in M}\mathrm{V}\mathrm{o}\mathrm{l}(B(x,1))>0.$
\end{thm}

      The necessity is a well known fact, see p.18 in~\cite{Hebey96}
 and was generalised by Carron~\cite{Carron94}.\ \par 
\quad The sufficiency of Theorem~\ref{SC18} was done by Varopoulos~\cite{Varopoulos89},
 see also Theorem 3.18, p. 37 in~\cite{Hebey96}, based on the
 work of Coulhon and Saloff-Coste~\cite{CoulSal93}.\ \par 
\ \par 
In this work we study the Sobolev embeddings for the case of
 vector bundles over $M.$\ \par 
\quad So let $G:=(H,\pi ,M)$ be a complex ${\mathcal{C}}^{m}$ vector
 bundle over $M$ of rank $N$ with fiber $H$ with a smooth scalar
 product $(\ ,\ ).$ We make the hypothesis that we have a \emph{metric}
 connection $\nabla ^{G}$ on $G,$ i.e. with $d$ the exterior
 derivation on $M,$ we have $d(u,v)=(\nabla ^{G}u,v)+(u,\nabla
 ^{G}v)$ for any smooth sections $u,\ v$ of $G.$ We shall call
 a vector bundle with these properties an \emph{adapted vector bundle.}\ \par 
\ \par 
\quad We shall introduce weights, given by the geometry, on a complete
 riemannian manifold $(M,g),$ in order to have Sobolev embeddings
 on vector bundles, \emph{always valid with these weights, without
 any curvature conditions nor volume control.}\ \par 
\ \par 
\quad In order to state the results, we need some definitions.\ \par 

\begin{defin}
~\label{SC27}Let $(M,g)$ be a riemannian manifold and $\displaystyle
 x\in M.$ We shall say that the geodesic ball $\displaystyle
 B(x,R)$ is $(0,\epsilon )$-{\bf admissible} if there is a chart
 $\displaystyle (B(x,R),\varphi )$ such that:\par 
(*)       $\displaystyle (1-\epsilon )\delta _{ij}\leq g_{ij}\leq
 (1+\epsilon )\delta _{ij}$ in $\displaystyle B(x,R)$ as bilinear form.\par 
We shall say that the geodesic ball $\displaystyle B(x,R)$ is
 $(1,\epsilon )$-{\bf admissible} if  moreover\par 
(**)     $\displaystyle \ \sum_{\left\vert{\beta }\right\vert
 =1}{R^{}\sup \ _{i,j=1,...,n,\ y\in B_{x}(R)}\left\vert{\partial
 ^{\beta }g_{ij}(y)}\right\vert }\leq \epsilon .$\par 
We shall denote ${\mathcal{A}}(0,\epsilon )$ the set of $(0,\epsilon
 )$-admissible balls on $M$ and  ${\mathcal{A}}(1,\epsilon )$
 the set of $(1,\epsilon )$-admissible balls on $M.$
\end{defin}
\quad The $(0,\epsilon )$-admissible balls will be adapted to functions
 and the $(1,\epsilon )$-admissible balls will be adapted to
 sections of the vector bundle $G.$\ \par 

\begin{defin}
~\label{SC28}Let $\displaystyle x\in M,$ we set $R'(x)=\sup \
 \lbrace R>0::B(x,R)\in {\mathcal{A}}(\epsilon )\rbrace .$ We
 shall say that $\displaystyle R_{\epsilon }(x):=\min \ (1,R'(x)/2)$
 is the $\epsilon $-{\bf admissible radius} at $\displaystyle x.$
\end{defin}
\quad Here  ${\mathcal{A}}(\epsilon )$ will be either ${\mathcal{A}}(0,\epsilon
 )$ or ${\mathcal{A}}(1,\epsilon )$ depending on the context,
 functions or sections of $G.$ The notation $R_{\epsilon }(x)$
 will means either $R_{0,\epsilon }(x)$ or $R_{1,\epsilon }(x)$
 depending on the choice of ${\mathcal{A}}(\epsilon ).$\ \par 
\quad Our general result with weights is:\ \par 

\begin{thm}
~\label{SB46}Let $(M,g)$ be a complete riemannian manifold of
 dimension $n.$ Let $G:=(H,\pi ,M)$ be a complex smooth adapted
 vector bundle over $M.$ Let $r\geq 1,\ m\geq k\geq 0$ and $\frac{1}{s}=\frac{1}{r}-\frac{(m-k)}{n}>0.$
 Let $\displaystyle w(x):=R_{\epsilon }(x)^{\gamma }$ and $\displaystyle
 w':=R_{\epsilon }(x)^{\nu }$ with  $\displaystyle \nu :=s(2+\gamma
 /r).$ Then $\displaystyle W_{G}^{m,r}(M,w)$ is embedded in $\displaystyle
 W_{G}^{k,s}(M,w')$ and:\par 
\vspace{.5em} \ \ \ \ \ \ \ \ \ \ \ \  $\displaystyle \forall
 u\in W_{G}^{m,r}(M,w),\ {\left\Vert{u}\right\Vert}_{W_{G}^{k,s}(M,w')}\leq
 C{\left\Vert{u}\right\Vert}_{W_{G}^{m,r}(M,w)}.$\vspace{.5em}  
\end{thm}
Because the admissible radius $R_{\epsilon }(x)$ is always smaller
 than $1,$ to get rid of the weights we just need that: $\exists
 \delta >0::\forall x\in M,\ R(x)\geq \delta .$\ \par 
\quad This is precisely what is given by a Theorem of Hebey-Herzlich
 ~\cite[Corollary, p. 7]{HebeyHerzlich97}, which has the easy Corollary:\ \par 

\begin{cor}
Let $(M,g)$ be a complete riemannian manifold. If the injectivity
 radius verifies $r_{inj}(x)\geq i>0$ and the Ricci curvature
 verifies $Rc_{(M,g)}(x)\geq \lambda g_{x}$ for some $\lambda
 \in {\mathbb{R}}$ and all $x\in M,$ then there exists a positive
 constant $\delta >0,$ depending only on $n,\epsilon ,\lambda
 ,i$ such that for any $x\in M,\ R_{0,\epsilon }(x)\geq \delta .$ \par 
\quad If moreover we have $\left\vert{Rc_{(M,g)}(x)}\right\vert \leq
 C$ for all $x\in M,$ then there exists a positive constant $\delta
 >0,$ depending only on $n,\epsilon ,i$ and  $C,$ such that for
 any $x\in M,\ R_{1,\epsilon }(x)\geq \delta .$
\end{cor}
\quad So the following results will be consequences of the Theorem~\ref{SB46}
 and of the Theorem of Hebey-Herzlich.\ \par 
\quad To state them precisely we shall need to weaken the definition
 of bounded geometry.\ \par 

\begin{defin}
A riemannian manifold $M$ has $k$-order {\bf weak bounded geometry}  if:\par 
\quad $\bullet $  the injectivity radius $r(x)$ at $x\in M$ is bounded
 below by some constant $\delta >0$ for any $\displaystyle x\in M$\par 
\quad $\bullet $ for $0\leq j\leq k,$ the covariant derivatives $\nabla
 ^{j}R_{c}$ of the Ricci curvature tensor are bounded in $L^{\infty }(M)$ norm.
\end{defin}
\quad Recall the classical Theorem of Cantor~\cite{Cantor74}:\ \par 

\begin{thm}
Let $(M,g)$ be a complete riemannian manifold. Let $G:=(H,\pi
 ,M)$ be a complex smooth adapted vector bundle over $M.$ Suppose
 $(M,g)$ verifies:\par 
\quad C1: the injectivity radius of $M$ is bounded away from zero.\par 
\quad C2: There is a $\delta $ such that for each $x\in M$ and $V,W\in
 T_{x}M,$ the
 sectional curvature $\left\vert{K_{x}(V,W)}\right\vert
 <\delta .$\par 
Let $\displaystyle 0\leq k<m$ and $\displaystyle 1/s=1/r-(m-k)/n.$
 Let $\displaystyle u\in W^{m,r}_{G}(M).$ Then we have $u\in
 W^{k,s}_{G}(M)$ with the control:\par 
\vspace{.5em} \ \ \ \ \ \ \ \ \ \ \ \  $\displaystyle {\left\Vert{u}\right\Vert}^{s}_{W^{k,s}_{G}(M)}\leq
 C{\left\Vert{u}\right\Vert}^{r}_{W^{m,r}_{G}(M)}.$\vspace{.5em}  
\end{thm}
Here we prove:\ \par 

\begin{thm}
Let $(M,g)$ be a complete riemannian manifold. Let $G:=(H,\pi
 ,M)$ be a complex smooth adapted vector bundle over $M.$ Suppose
 $(M,g)$ has a $0$-order weak bounded geometry. Let $\displaystyle
 0\leq k<m$ and $\displaystyle 1/s=1/r-(m-k)/n.$ Let $\displaystyle
 u\in W^{m,r}_{G}(M).$ Then we have $u\in W^{k,s}_{G}(M)$ with the control:\par 
\vspace{.5em} \ \ \ \ \ \ \ \ \ \ \ \  $\displaystyle {\left\Vert{u}\right\Vert}^{s}_{W^{k,s}_{G}(M)}\leq
 C{\left\Vert{u}\right\Vert}^{r}_{W^{m,r}_{G}(M)}.$\vspace{.5em}  \par 
In the case of functions instead of sections of $G,$ the conditions
 are weaker: if the injectivity radius of $M$ is bounded away
 from zero and if the Ricci curvature verifies $Rc_{(M,g)}(x)\geq
 \lambda g_{x}$ for some $\lambda \in {\mathbb{R}}$ and all $x\in
 M,$ we get\par 
\vspace{.5em} \ \ \ \ \ \ \ \ \ \ \ \  $\displaystyle {\left\Vert{u}\right\Vert}^{s}_{W^{k,s}(M)}\leq
 C{\left\Vert{u}\right\Vert}^{r}_{W^{m,r}(M)}.$\vspace{.5em}  
\end{thm}
\quad This improves the Theorem by Cantor because he used the hypothesis
 that all the sectional curvatures are bounded and here we need
 only that the Ricci curvature be bounded.\ \par 
\ \par 
\quad We also prove a global Gaffney's type inequality:\ \par 

\begin{thm}
Let $(M,g)$ be a complete riemannian manifold. Let $d$ be the
 exterior derivation on $M$ and $d^{*}$ its formal adjoint. Let
 $\omega $ be a $p$-differential form in $M.$ If $(M,g)$ has
 a $0$-order weak bounded geometry, then we have, with $r\geq 1$:\par 
\vspace{.5em} \ \ \ \ \ \ \ \ \ \ \ \  $\displaystyle {\left\Vert{\omega
 }\right\Vert}_{W_{p}^{1,r}(M)}\leq C({\left\Vert{d(\omega )}\right\Vert}_{L_{p+1}^{r}(M)}+{\left\Vert{d^{*}(\omega
 )}\right\Vert}_{L_{p-1}^{r}(M)}+{\left\Vert{\omega }\right\Vert}_{L_{p}^{r}(M)}).$\vspace{.5em}
  
\end{thm}
And, as an easy corollary:\ \par 

\begin{cor}
Let $(M,g)$ be a complete riemannian manifold with a $0$-order
 weak bounded geometry.  Let $r\geq 1$and let $\omega $ be a
 $p$-differential form in $M.$ We have:\par 
\vspace{.5em} \ \ \ \ \ \ \ \ \ \ \ \  $\displaystyle {\left\Vert{\omega
 }\right\Vert}_{L_{p}^{s}(M)}\leq C({\left\Vert{d(\omega )}\right\Vert}_{L_{p+1}^{r}(M)}+{\left\Vert{d^{*}(\omega
 )}\right\Vert}_{L_{p-1}^{r}(M)}+{\left\Vert{\omega }\right\Vert}_{L_{p}^{r}(M)})$\vspace{.5em}
  \par 
with $\frac{1}{s}=\frac{1}{r}-\frac{1}{n}>0.$
\end{cor}
\quad N. Lohou\'e~\cite{Lohoue85} proved the same result under the
 stronger hypothesis that $(M,g)$ has a $2$-order bounded geometry
 plus some other hypotheses on the laplacian and the range of
 $r.$ Here already $0$-order weak bounded geometry is enough.\ \par 
\ \par 
\quad This work is presented as follow:\ \par 
$\bullet $ In the next Section we study the main properties of
 the $\epsilon $-admissible balls.\ \par 
$\bullet $ In Section~\ref{SB37} we define the vector bundle
 $G$ we are interested in and precise the metric connexion $\nabla
 ^{G}$ on it we shall use.\ \par 
$\bullet $ In Section~\ref{SB45} we define the Sobolev spaces
 of smooth sections of $G,$ with weights. We prove in this Section
 a generalisation of a nice result of T. Aubin~\cite{Aubin82}
 which says that, in order to prove Sobolev embeddings for sections
 of $G$ with weights, we have just to prove them at the first
 level. This is crucial for our estimates.\ \par 
$\bullet $ In Section~\ref{SB38} we prove the local estimates,
 i.e. for sections of $G$ in the $\epsilon $-admissible balls.\ \par 
$\bullet $ In SubSection~\ref{SB39} we also prove a local Gaffney
 type inequality, using a result by C. Scott~\cite{Scott95}.\ \par 
$\bullet ${\it  }Then in order to get global results we group
 the $\epsilon $-admissible balls via a Vitali type covering
 in Section~\ref{SB36}.\ \par 
$\bullet $ In Section~\ref{SB40} we prove the global estimates
 for functions, sections of $G$ and Gaffney type in $L^{r}.$\ \par 
$\bullet $ In Section~\ref{SB41} we improve the classical Sobolev
 embeddings to the case of riemannian manifolds with weak bounded
 geometry, by use of a Theorem of Hebey-Herzlich~\cite{HebeyHerzlich97}.
 This implies the validity of Sobolev embeddings for vector bundles
 in compact riemannian manifold without boundary.\ \par 
$\bullet $ In SubSection~\ref{SB42} we deduce from the compact
 case without boundary the validity of Sobolev embeddings for
 vector bundles in compact riemannian manifold with smooth boundary.
 We use here the method of the "double" manifold.\ \par 
$\bullet $ Finally in SubSection~\ref{SB43} we introduce the
 \emph{lifted doubling property} and we study the case of hyperbolic
 manifolds.\ \par 

\section{Admissible balls.~\label{S11}}
\quad Recall the definition of the admissible radius:\ \par 

\begin{defin}
Let $\displaystyle x\in M,$ we set $R'(x)=\sup \ \lbrace R>0::B(x,R)\in
 {\mathcal{A}}(\epsilon )\rbrace .$ We shall say that $\displaystyle
 R_{\epsilon }(x):=\min \ (1,R'(x)/2)$ is the $\epsilon $-{\bf
 admissible radius} at $\displaystyle x.$
\end{defin}
Then clearly if $B(x,R)\in {\mathcal{A}}(\epsilon ),$ i.e. is
 $\epsilon $-admissible, then so is $B(x,S)$ if $S\leq R.$ \ \par 

\begin{rem}
Let $\displaystyle x,y\in M.$ Suppose that $R'(x)>d_{g}(x,y),$
 where $d_{g}(x,y)$ is the riemannian distance between $x$ and
 $y.$ Consider the ball $B(y,\rho )$ of center $y$ and radius
 $\rho :=R'(x)-d_{g}(x,y).$ This ball is contained in $B(x,R'(x))$
 hence, by definition of $R'(x),$ we have that all the points
 in $B(y,\rho )$ verify the conditions 1) and 2) so, by definition
 of $R'(y),$ we have that $R'(y)\geq R'(x)-d_{g}(x,y).$ If $\displaystyle
 R'(x)\leq d_{g}(x,y)$ this is also true because $\displaystyle
 R'(y)>0.$ Exchanging $x$ and $y$ we get that $\left\vert{R'(y)-R'(x)}\right\vert
 \leq d_{g}(x,y).$\par 
\quad Hence $R'(x)$ is $1$-lipschitzian so it is continuous. So the
 $\epsilon $-admissible radius $R_{\epsilon }(x)$ is also continuous.
\end{rem}

\begin{rem}
~\label{SC16} Because on our admissible ball $B(x,R_{\epsilon
 }(x))$ there is a diffeomorphism from $\displaystyle B(x,R)$
 to $\varphi (B(x,R))\subset {\mathbb{R}}^{n},$ i.e. on an open
 set in the tangent space $T_{x}M,$ we get that the injectivity
 radius $r_{inj}(x)$ always verifies $\displaystyle r_{inj}(x)\geq
 R_{\epsilon }(x).$
\end{rem}
\ \par 

\begin{lem}
~\label{S2}(Slow variation of the admissible radius) Let $(M,g)$
 be a riemannian manifold then with $R(x)=R_{\epsilon }(x)=$
 the $\epsilon $-admissible radius at $x\in M.$ We get:\par 
\quad \quad $\forall y\in B(x,R(x))$ we have $R(x)/2\leq R(y)\leq 2R(x).$
\end{lem}
\quad Proof.\ \par 
Let $x,y\in M$ and $d(x,y)$ the riemannian distance on $(M,g).$
 Let $y\in B(x,R(x))$ then $d(x,y)\leq R(x)$ and suppose first
 that $R(x)\geq R(y).$ Then, because $R(x)=R'(x)/2,$ we get $y\in
 B(x,R'(x)/2)$ hence we have $B(y,R'(x)/2)\subset B(x,R'(x)).$
 But by the definition of $R'(x),$ the ball $\displaystyle B(x,R'(x))$
 is admissible and this implies that the ball $\displaystyle
 B(y,R'(x)/2)$ is also admissible for exactly the same constants
 and the same chart; this implies that $R'(y)\geq R'(x)/2$ hence
 $R(y)\geq R(x)/2,$ so $R(x)\geq R(y)\geq R(x)/2.$\ \par 
\quad If $R(x)\leq R(y)$ then $d(x,y)\leq R(x)\Rightarrow d(x,y)\leq
 R(y)\Rightarrow x\in B(y,R'(y)/2)\Rightarrow B(x,R'(y)/2)\subset
 B(y,R'(y)).$ Hence the same way as above we get $R(y)\geq R(x)\geq
 R(y)/2\Rightarrow R(y)\leq 2R(x).$ So in any case we proved that\ \par 
\quad \quad \quad $\forall y\in B(x,R(x))$ we have $R(x)/2\leq R(y)\leq 2R(x).$
 $\hfill\blacksquare $\ \par 

\section{Vector bundle.~\label{SB37}}
\quad Let $(M,g)$ be a complete riemannian manifold and let $G:=(H,\pi
 ,M)$ be an adapted complex ${\mathcal{C}}^{m}$ vector bundle
 over $M$ of rank $N$ with fiber $H.$ Recall that this means
 that $G$ has a smooth scalar product $(\ ,\ )$ and  a \emph{metric}
 connection $\nabla ^{G}:{\mathcal{C}}^{\infty }(M,G)\rightarrow
 {\mathcal{C}}^{\infty }(M,G\otimes T^{*}M),$ i.e. verifying
 $\displaystyle d(u,v)=(\nabla ^{G}u,v)+(u,\nabla ^{G}v),$ where
 $d$ is the exterior derivative on $M.$ See ~\cite[Section 13]{TaylorGD}.\ \par 

\begin{lem}
~\label{SB30} The $\epsilon $-admissible balls $B(x,R_{\epsilon
 }(x))$ trivialise the bundle $G.$
\end{lem}
\quad Proof.\ \par 
Because if $B(x,R)$ is an $\epsilon $-admissible ball, we have
 by Remark~\ref{SC16} that $R\leq r_{inj}(x).$ Then, one can
 choose a local frame field for $G$ on $\displaystyle B(x,R)$
 by radial parallel translation, as done in~\cite[Section 13,
 p.86-87]{TaylorGD}, see also~\cite[p. 4, eq. (1.3)]{Nistor06}.
 This means that the $\epsilon $-admissible ball also trivialises
 the bundle $G.$ $\hfill\blacksquare $\ \par 
\ \par 
\quad If $\partial _{j}:=\partial /\partial x_{j}$ in a coordinate
 system on, say $B(x_{0},R),$ and with a local frame $\lbrace
 e_{\alpha }\rbrace _{\alpha =1,...,N},$ we have, for a smooth
 sections of $G,$ $u=u^{\alpha }e_{\alpha }$ with the Einstein
 summation convention. We set:\ \par 
\vspace{.5em} \ \ \ \ \ \ \ \ \ \ \ \  $\displaystyle \ \nabla
 _{\partial _{j}}u=(\partial _{j}u^{\alpha }+u^{\beta }\Gamma
 ^{G,\alpha }_{\beta j})e_{\alpha },$\vspace{.5em}  \ \par 
the Christoffel coefficients $\displaystyle \Gamma ^{G,\alpha
 }_{\beta j}$ being defined by $\displaystyle \nabla _{\partial
 _{j}}e_{\beta }=\Gamma ^{G,\alpha }_{\beta j}e_{\alpha }.$\ \par 
\quad We shall make the following Control by the Metric Tensor hypothesis,
 for $B(x_{0},R)\in {\mathcal{A}}(1,\epsilon )$:\ \par 
(CMT)        $\forall x\in B(x_{0},R),\ \left\vert{\Gamma ^{G,\alpha
 }_{\beta j}(x)}\right\vert \leq C(n,G,\epsilon )\sum_{\left\vert{\beta
 }\right\vert =1}{\sup \ _{i,j=1,...,n,}\left\vert{\partial ^{\beta
 }g_{ij}(x)}\right\vert },$\ \par 
the constant $C$ depending only on $n,\epsilon $ and $G.$\ \par 
\quad This hypothesis is natural in the sense that it is true for tensor
 bundles over $M.$\ \par 

\begin{lem}
Let $F$ be a tensor bundle over $M.$ Then the hypothesis (CMT) is true.
\end{lem}
\quad Proof.\ \par 
Let $\displaystyle \Gamma ^{k}_{lj}$ be the Christoffel coefficients
 of the Levi-Civita connexion on the tangent bundle $TM.$ We have\ \par 
\quad \quad \quad \begin{equation} \Gamma ^{k}_{ij}=\frac{1}{2}g^{il}(\frac{\partial
 g_{kl}}{\partial x^{j}}+\frac{\partial g_{lj}}{\partial x^{k}}-\frac{\partial
 g_{jk}}{\partial x^{l}}).\label{SB34}\end{equation}\ \par 
Now on $B(x_{0},R)\in {\mathcal{A}}(\epsilon ),$ we have $(1-\epsilon
 )\delta _{ij}\leq g_{ij}\leq (1+\epsilon )\delta _{ij}$ as bilinear
 forms. Hence\ \par 
\vspace{.5em} \ \ \ \ \ \ \ \ \ \ \ \  $\displaystyle \forall
 x\in B(x_{0},R),\ \ \left\vert{\Gamma ^{k}_{ij}(x)}\right\vert
 \leq \frac{3}{2}(1+\epsilon )\sum_{\left\vert{\beta }\right\vert
 =1}{\sup \ _{i,j=1,...,n,}\left\vert{\partial ^{\beta }g_{ij}(x)}\right\vert
 }$\vspace{.5em}  \ \par 
in a coordinates chart on $B(x_{0},R).$\ \par 
\quad Let $\omega :=\omega _{k}dx^{k}$ be a $1$-form in this chart.
 Then $\nabla _{\partial _{j}}\omega =(\partial _{j}\omega _{k}-\Gamma
 ^{l}_{kj}a_{l})dx^{k}.$ In particular we have $\nabla _{\partial
 _{j}}dx^{m}=(-\Gamma ^{m}_{kj})dx^{k}.$ Hence (CMT) is true.
 The same if $X$ is a vector field, $X:=X^{k}\partial _{k}$ we
 have $\nabla _{\partial _{j}}X=(\partial _{j}X^{k}+\Gamma ^{k}_{lj}X^{l})\partial
 _{k}.$ So $\nabla _{\partial _{j}}\partial _{l}=(\Gamma ^{k}_{lj})\partial
 _{k}$ and (CMT) is still true.\ \par 
\quad Now on $(p,q)$ tensors, using the fact that $\nabla $ is a derivation,
 we have that\ \par 
\vspace{.5em} \ \ \ \ \ \ \ \ \ \ \ \  $\displaystyle \nabla
 _{\partial _{j}}(dx^{k_{1}}\otimes \cdot \cdot \cdot \otimes
 dx^{k_{q}}\otimes \partial _{l_{1}}\otimes \cdot \cdot \cdot
 \otimes \partial _{l_{p}})=$\vspace{.5em}  \ \par 
\vspace{.5em} \ \ \ \ \ \ \ \ \ \ \ \ \ \ \ \ \ \ \ \ \ \ \ \
  $\displaystyle =(\nabla _{\partial _{j}}dx^{k_{1}})\otimes
 \cdot \cdot \cdot \otimes dx^{k_{q}}\otimes \partial _{l_{1}}\otimes
 \cdot \cdot \cdot \otimes \partial _{l_{p}}+\cdot \cdot \cdot
 $\vspace{.5em}  \ \par 
\vspace{.5em} \ \ \ \ \ \ \ \ \ \ \ \ \ \ \ \ \ \ \ \ \ \ \ \
 \ \ \ \ \ \  $\displaystyle \cdot \cdot \cdot +dx^{k_{1}}\otimes
 \cdot \cdot \cdot \otimes dx^{k_{q}}\otimes \partial _{l_{1}}\otimes
 \cdot \cdot \cdot \otimes (\nabla _{\partial _{j}}\partial _{l_{p}})$\vspace{.5em}
  \ \par 
is a linear combination of $\displaystyle \Gamma ^{k}_{lj}$ with
 constant coefficients, hence the (CMT) is also true, with a
 constant depending only on the dimension $n$ of $M$ and $p,q,\epsilon .$\ \par 
\quad A special mention for the bundle $\Lambda ^{p}(M)$ of $p$-forms
 on $M,$ which is a sub bundle of the $(0,p)$ tensors bundle.
 It also has the (CMT) property.\ \par 
The proof is complete. $\hfill\blacksquare $ \ \par 

\section{Sobolev spaces for sections of $G$ with weight.~\label{SB45}}
\quad We have seen that $\nabla ^{G}:\ {\mathcal{C}}^{\infty }(M,G)\rightarrow
 {\mathcal{C}}^{\infty }(M,G\otimes T^{*}M).$ On the tensor product
 of two Hilbert spaces we put the canonical scalar product $(u\otimes
 \omega ,v\otimes \mu ):=(u,v)(\omega ,\mu ),$ with $u\otimes
 \omega \in G\otimes T^{*}M,$ and completed by linearity to all
 elements of the tensor product. On $T^{*}M$ we have the Levi-Civita
 connexion $\nabla ^{M},$ which is of course a metric one, and
 on $G$ we have the metric connexion $\nabla ^{G}$ so we define
 a connexion on the tensor product $G\otimes T^{*}M$:\ \par 
\vspace{.5em} \ \ \ \ \ \ \ \ \ \ \ \  $\displaystyle \nabla
 ^{G\otimes T^{*}M}(u\otimes \omega )=(\nabla ^{G}u)\otimes \omega
 +u\otimes (\nabla ^{T^{*}M}\omega )$\vspace{.5em}  \ \par 
by asking that this connexion be a derivation. We get easily that\ \par 
\vspace{.5em} \ \ \ \ \ \ \ \ \ \ \ \  $\displaystyle \nabla
 ^{G\otimes T^{*}M}\ :\ {\mathcal{C}}^{\infty }(M,G\otimes T^{*}M)\rightarrow
 {\mathcal{C}}^{\infty }(M,G\otimes (T^{*}M)^{\otimes 2})$\vspace{.5em}  \ \par 
is still a \emph{metric} connexion, i.e.\ \par 
\vspace{.5em} \ \ \ \ \ \ \ \ \ \ \ \  $\displaystyle \ d(u\otimes
 \omega ,v\otimes \mu )=(\nabla ^{G\otimes T^{*}M}(u\otimes \omega
 ),\ v\otimes \mu )+(u\otimes \omega ,\ \nabla ^{G\otimes T^{*}M}(v\otimes
 \mu )).$\vspace{.5em}  \ \par 
\quad We define by iteration $\nabla ^{j}u:=\nabla (\nabla ^{j-1}u)$
 on the section $u$ of $G$ and the associated pointwise scalar
 product $(\nabla ^{j}u(x),\nabla ^{j}v(x))$ which is defined
 on $G\otimes (T^{*}M)^{\otimes j},$ with again the metric connection\ \par 
\vspace{.5em} \ \ \ \ \ \ \ \ \ \ \ \  $\displaystyle d(\nabla
 ^{j}u,\nabla ^{j}v)(x)=(\nabla ^{j+1}u,\nabla ^{j}v)(x)+(\nabla
 ^{j}u,\nabla ^{j+1}v)(x).$\vspace{.5em}  \ \par 
\ \par 
\quad Let $w$ be a weight on $M,$ i.e. a positive measurable function
 on $M.$ If $k\in {\mathbb{N}}$ and $r\geq 1$ are given, we denote
 by ${\mathcal{C}}^{k,r}_{G}(M,w)$ the space of smooth sections
 of $G\ \omega \in {\mathcal{C}}^{\infty }(M)$ such that $\ \left\vert{\nabla
 ^{j}\omega }\right\vert \in L^{r}(M,w)$ for $j=0,...,k$ with
 the pointwise modulus associated to the pointwise scalar product. Hence\ \par 
\vspace{.5em} \ \ \ \ \ \ \ \ \ \ \ \  $\displaystyle {\mathcal{C}}^{k,r}_{G}(M,w):=\lbrace
 \omega \in {\mathcal{C}}_{G}^{\infty }(M),\ \forall j=0,...,k,\
 \int_{M}{\left\vert{\nabla ^{j}\omega }\right\vert ^{r}(x)w(x)dv(x)}<\infty
 \rbrace ,$\vspace{.5em}  \ \par 
with $dv$ the volume measure on $(M,g).$\ \par 
\quad Now we have, see M. Cantor ~\cite[Definition 1 \&  2, p. 240]{Cantor74}
 for the case without weight:\ \par 

\begin{defin}
The Sobolev space $W_{G}^{k,r}(M,w)$ is the completion of ${\mathcal{C}}^{k,r}_{G}(M,w)$
 with respect to the norm:\par 
\vspace{.5em} \ \ \ \ \ \ \ \ \ \ \ \  $\displaystyle {\left\Vert{\omega
 }\right\Vert}_{W_{G}^{k,r}(M,w)}=\sum_{j=0}^{k}{{\left({\int_{M}{\left\vert{\nabla
 ^{j}\omega (x)}\right\vert ^{r}w(x)dv(x)}}\right)}^{1/r}}.$\vspace{.5em}  
\end{defin}
The usual case is when $w\equiv 1.$ Then we write simply $W_{G}^{k,r}(M).$\
 \par 
\ \par 
\quad We shall apply these well known facts to generalise a nice result
 of T. Aubin.\ \par 
Let $w(x),\ w'(x)$ be  weights on the complete riemannian manifold
 $(M,g).$ We have:\ \par 

\begin{prop}
~\label{SC14}Let $(M,g)$ be a complete riemannian manifold. If
 $W_{G}^{1,r_{0}}(M,w)$ is embedded in $L_{G}^{s_{0}}(M,w'),$
 with $\frac{1}{s_{0}}=\frac{1}{r_{0}}-\frac{1}{n}\ (1\leq r_{0}<n),$
 then $W_{G}^{k,r}(M,w)$ is embedded in $W_{G}^{l,s_{l}}(M,w'),$
 with $\frac{1}{s_{l}}=\frac{1}{r}-\frac{(k-l)}{n}>0.$
\end{prop}
\quad Proof.\ \par 
We shall copy the proof of Proposition 2.11, p. 36 in~\cite{Aubin82},
 replacing $L^{s}(M)$ by $L_{G}^{s}(M,w')$ and $W^{k,r}(M)$ by
 $W_{G}^{k,r}(M,w)$ and simplifying a little bit the argument.\ \par 
\quad Let $m$ be an integer and let $\omega \in {\mathcal{C}}_{G}^{m+1}.$
 We have pointwise:\ \par 
\quad \quad \quad \begin{equation} \ \left\vert{\nabla \left\vert{\nabla ^{m}\omega
 }\right\vert }\right\vert \leq \left\vert{\nabla ^{m+1}\omega
 }\right\vert .\label{SC12}\end{equation}\ \par 
To see this, by $\left\vert{\nabla ^{m}\omega }\right\vert ^{2}(x)=(\nabla
 ^{m}\omega ,\nabla ^{m}\omega )(x),$ we have\ \par 
\vspace{.5em} \ \ \ \ \ \ \ \ \ \ \ \  $\displaystyle \nabla
 \left\vert{\nabla ^{m}\psi }\right\vert ^{2}=\nabla (\nabla
 ^{m}\omega ,\nabla ^{m}\omega )(x)=2(\nabla ^{m+1}\omega ,\nabla
 ^{m}\omega ),$\vspace{.5em}  \ \par 
the last equality because $\nabla $ is a metric connection.\ \par 
By the Cauchy-Schwartz inequality, we get\ \par 
\vspace{.5em} \ \ \ \ \ \ \ \ \ \ \ \  $\displaystyle \ \left\vert{\nabla
 \left\vert{\nabla ^{m}\omega }\right\vert ^{2}}\right\vert (x)\leq
 2\left\vert{\nabla ^{m+1}\omega }\right\vert \left\vert{\nabla
 ^{m}\omega }\right\vert (x).$\vspace{.5em}  \ \par 
Setting $F(x):=\left\vert{\nabla ^{m}\psi }\right\vert (x),$
 because $\nabla $ is a derivation, we also have $\nabla F^{2}(x)=2F(x)\nabla
 F(x),$ hence:\ \par 
\vspace{.5em} \ \ \ \ \ \ \ \ \ \ \ \  $\displaystyle \left\vert{\nabla
 \left\vert{\nabla ^{r}\psi }\right\vert ^{2}}\right\vert =2\left\vert{F}\right\vert
 \left\vert{\nabla \left\vert{F}\right\vert }\right\vert \leq
 2\left\vert{\nabla ^{m+1}\omega }\right\vert \left\vert{\nabla
 ^{m}\omega }\right\vert ,$\vspace{.5em}  \ \par 
so $\left\vert{\nabla \left\vert{\nabla ^{m}\omega }\right\vert
 }\right\vert \leq \left\vert{\nabla ^{m+1}\omega }\right\vert .$\ \par 
\ \par 
\quad Since $W_{G}^{1,r_{0}}(M,w)$  is embedded in $L_{G}^{s_{0}}(M,w'),$
  there exists a constant $A,$ such that for all $\varphi \in
 W_{G}^{1,r_{0}}(M,w)$: (for now on, we do not indicate the subscript
 to ease the notation.)\ \par 
\vspace{.5em} \ \ \ \ \ \ \ \ \ \ \ \  $\displaystyle \ {\left\Vert{\varphi
 }\right\Vert}_{L^{s_{0}}(M,w')}\leq A({\left\Vert{\nabla \varphi
 }\right\Vert}_{L^{r_{0}}(M,w)}+{\left\Vert{\varphi }\right\Vert}_{L^{r_{0}}(M,w)}).$\vspace{.5em}
  \ \par 
Let us apply this inequality with $\varphi =\left\vert{\nabla
 ^{m}\omega }\right\vert ,$ assuming $\varphi ,$ which is a function
 now, belongs to $W^{1,r_{0}}(M,w)$:\ \par 
\vspace{.5em} \ \ \ \ \ \ \ \ \ \ \ \  $\displaystyle {\left\Vert{\nabla
 ^{m}\omega }\right\Vert}_{L^{s_{0}}(M,w')}\leq A({\left\Vert{\nabla
 \left\vert{\nabla ^{m}\omega }\right\vert }\right\Vert}_{L^{r_{0}}(M,w)}+{\left\Vert{\nabla
 ^{m}\omega }\right\Vert}_{L^{r_{0}}(M,w)})$\vspace{.5em}  \ \par 
hence, using~(\ref{SC12}) and integrating,\ \par 
\quad \quad \quad \begin{equation} {\left\Vert{\nabla ^{m}\omega }\right\Vert}_{L^{s_{0}}(M,w')}\leq
 A({\left\Vert{\nabla ^{m+1}\omega }\right\Vert}_{L^{r_{0}}(M,w)}+{\left\Vert{\nabla
 ^{m}\omega }\right\Vert}_{L^{r_{0}}(M,w)}).\label{SC13}\end{equation}\ \par 
Now let $\omega \in W_{G}^{k,r}(M,w)\cap {\mathcal{C}}_{G}^{\infty
 }(M).$ Applying inequalities~(\ref{SC12}) and~(\ref{SC13}) with
 $r=r_{0}$ and $m=k-1,\ k-2,...,$ we find, for $\nabla ^{j}\omega ,$\ \par 
\vspace{.5em} \ \ \ \ \ \ \ \ \ \ \ \  $\displaystyle \ {\left\Vert{\nabla
 ^{k-1}\omega }\right\Vert}_{L^{s_{k-1}}(M,w')}\leq A({\left\Vert{\nabla
 ^{k}\omega }\right\Vert}_{L^{r}(M,w)}+{\left\Vert{\nabla ^{k-1}\omega
 }\right\Vert}_{L^{r}(M,w)}),$\vspace{.5em}  \ \par 
\vspace{.5em} \ \ \ \ \ \ \ \ \ \ \ \  $\displaystyle {\left\Vert{\nabla
 ^{k-2}\omega }\right\Vert}_{L^{s_{k-1}}(M,w')}\leq A({\left\Vert{\nabla
 ^{k-1}\omega }\right\Vert}_{L^{r}(M,w)}+{\left\Vert{\nabla ^{k-2}\omega
 }\right\Vert}_{L^{r}(M,w)}),$\vspace{.5em}  \ \par 
\vspace{.5em} \ \ \ \ \ \ \ \ \ \ \ \  $\displaystyle {\left\Vert{\omega
 }\right\Vert}_{L^{s_{k-1}}(M,w')}\leq A({\left\Vert{\nabla \omega
 }\right\Vert}_{L^{r}(M,w)}+{\left\Vert{\psi }\right\Vert}_{L^{r}(M,w)}).$\vspace{.5em}
  \ \par 
Thus\ \par 
\vspace{.5em} \ \ \ \ \ \ \ \ \ \ \ \  $\displaystyle \ {\left\Vert{\omega
 }\right\Vert}_{W^{k-1,s_{k-1}}(M,w')}\leq 2A{\left\Vert{\omega
 }\right\Vert}_{W^{k,r}(M,w)}.$\vspace{.5em}  \ \par 
Therefore a Cauchy sequence in $W_{G}^{k,r}(M,w)$ of sections
 of $G$  is a Cauchy sequence in $W_{G}^{k-1,s_{k-1}}(M,w'),$
 and the preceding inequality holds for all $\psi \in W_{G}^{k,r}(M,w)$
 and we get:\ \par 
\vspace{.5em} \ \ \ \ \ \ \ \ \ \ \ \  $\displaystyle W_{G}^{k,r}(M,w)\subset
 W_{G}^{k-1,s_{k-1}}(M,w').$\vspace{.5em}  \ \par 
Now with $w=w'$ we prove similarly the following embeddings:\ \par 
\vspace{.5em} \ \ \ \ \ \ \ \ \ \ \ \  $\displaystyle W_{G}^{k-1,s_{k-1}}(M,w')\subset
 W_{G}^{k-2,s_{k-2}}(M,w')\subset \cdot \cdot \cdot \subset W_{G}^{l,s_{l}}(M,w').$\vspace{.5em}
  \ \par 
Hence $W_{G}^{k,r}(M,w)\subset W_{G}^{l,s_{l}}(M,w').$ $\hfill\blacksquare
 $\ \par 
\ \par 
\quad Proposition~\ref{SC14} says that, in order to prove Sobolev embeddings
 with weights, we have just to prove that $W_{G}^{1,r}(M,w)$
  is embedded in $L_{G}^{s}(M,w'),$ with $\frac{1}{s}=\frac{1}{r_{}}-\frac{1}{n}\
 (1\leq r<n).$\ \par 
\quad This is very important here because we shall have just to deal
 with first order Sobolev spaces. Hence we have to work only
 with $\nabla ^{G}u$ which, by our assumption (CMT), implies
 at most the first order derivatives of the metric tensor.\ \par 
\ \par 
\quad The aim now is to prove that $W_{G}^{1,r_{0}}(M,w)$  is embedded
 in $L_{G}^{s_{0}}(M,w'),$ with $\frac{1}{s_{0}}=\frac{1}{r_{0}}-\frac{1}{n}\
 (1\leq r_{0}<n)$ then we shall be able to apply Proposition~\ref{SC14}.\ \par 

\section{Local estimates.~\label{SB38}}

\subsection{Sobolev comparison estimates for functions.}

\begin{lem}
~\label{HCS43}We have the Sobolev comparison estimates where
 $B(x,R)$ is a $(0,\epsilon )$-admissible ball in $M$ and $\varphi
 \ :\ B(x,R)\rightarrow {\mathbb{R}}^{n}$  is the admissible
 chart relative to $B(x,R),$\par 
\vspace{.5em} \ \ \ \ \ \ \ \ \ \ \ \  $\displaystyle \forall
 u\in W^{1,r}(B(x,R)),\ {\left\Vert{u}\right\Vert}_{W^{1,r}(B(x,R))}\leq
 C{\left\Vert{u\circ \varphi ^{-1}}\right\Vert}_{W^{1,r}(\varphi
 (B(x,R)))},$\vspace{.5em}  \par 
and, with $B_{e}(0,t)$ the euclidean ball in ${\mathbb{R}}^{n}$
 centered at $0$ and of radius $t,$\par 
\vspace{.5em} \ \ \ \ \ \ \ \ \ \ \ \  $\displaystyle {\left\Vert{v}\right\Vert}_{W^{1,r}(B_{e}(0,(1-\epsilon
 )R))}\leq C'{\left\Vert{u}\right\Vert}_{W^{1,r}(B(x,R))}.$\vspace{.5em}  \par 
The constants $C,C'$ depending only on $\epsilon ,n$ and not
 on $B\in {\mathcal{A}}(\epsilon ).$
\end{lem}
\quad Proof.\ \par 
We have to compare the norms of $u,\ \nabla u,$ with the corresponding
 ones for $v:=u\circ \varphi ^{-1}$ in ${\mathbb{R}}^{n}.$\ \par 
First we have because $(1-\epsilon )\delta _{ij}\leq g_{ij}\leq
 (1+\epsilon )\delta _{ij}$ in $B(x,R)$:\ \par 
\vspace{.5em} \ \ \ \ \ \ \ \ \ \ \ \  $\displaystyle B_{e}(0,(1-\epsilon
 )R)\subset \varphi (B(x,R))\subset B_{e}(0,(1+\epsilon )R).$\vspace{.5em}
  \ \par 
For functions, in a coordinates chart, we have $(\nabla u)_{j}:=\partial
 _{j}u$ hence with $\forall y\in B(x,R),\ z:=\varphi (y),$ we get:\ \par 
\vspace{.5em} \ \ \ \ \ \ \ \ \ \ \ \  $\displaystyle \forall
 y\in B(x,R),\ \left\vert{u(y)}\right\vert =\left\vert{v(z)}\right\vert
 ,\ \ \left\vert{\nabla u(y)}\right\vert \leq \left\vert{\partial
 v(z)}\right\vert .$\vspace{.5em}  \ \par 
\quad Integrating this we get:\ \par 
\vspace{.5em} \ \ \ \ \ \ \ \ \ \ \ \  $\displaystyle {\left\Vert{\nabla
 u}\right\Vert}_{L^{r}(B(x,R))}\leq {\left\Vert{\partial v}\right\Vert}_{L^{r}(B_{e}(0,(1+\epsilon
 )R))}.$\vspace{.5em}  \ \par 
So, using that\ \par 
\vspace{.5em} \ \ \ \ \ \ \ \ \ \ \ \  $\displaystyle {\left\Vert{u}\right\Vert}_{W^{1,r}(B(x,R))}={\left\Vert{\nabla
 u}\right\Vert}_{L^{r}(B(x,R))}+{\left\Vert{u}\right\Vert}_{L^{r}(B(x,R))},$\vspace{.5em}
  \ \par 
we get\ \par 
\vspace{.5em} \ \ \ \ \ \ \ \ \ \ \ \  $\displaystyle {\left\Vert{u}\right\Vert}_{W^{1,r}(B(x,R))}\leq
 C({\left\Vert{\partial v}\right\Vert}_{L^{r}(B_{e}(0,(1+\epsilon
 )R))}+{\left\Vert{v}\right\Vert}_{L^{r}(B_{e}(0,(1+\epsilon
 )R))}),$\vspace{.5em}  \ \par 
hence\ \par 
\vspace{.5em} \ \ \ \ \ \ \ \ \ \ \ \  $\displaystyle {\left\Vert{u}\right\Vert}_{W^{1,r}(B(x,R))}\leq
 C{\left\Vert{v}\right\Vert}_{W^{1,r}(\varphi (B(x,R))}.$\vspace{.5em}  \ \par 
Of course all these estimates can be reversed so we also have\ \par 
\vspace{.5em} \ \ \ \ \ \ \ \ \ \ \ \  $\displaystyle {\left\Vert{v}\right\Vert}_{W^{1,r}(B_{e}(0,(1-\epsilon
 )R))}\leq C{\left\Vert{u}\right\Vert}_{W^{1,r}(B(x,R))}.$\vspace{.5em}  \ \par 
\quad This ends the proof of the lemma. $\hfill\blacksquare $\ \par 
\ \par 
\quad We have to study the behavior of the Sobolev embeddings w.r.t.
 the radius. Set $B_{R}:=B_{e}(0,R)$ an euclidean ball in ${\mathbb{R}}^{n}.$
 For this purpose we have by ~\cite[Lemma 7.7]{ellipticEq18}
 in the special case $m=1.$
\begin{lem}
~\label{3S3}We have, with $t::\frac{1}{t}=\frac{1}{r}-\frac{1}{n},$\par 
\vspace{.5em} \ \ \ \ \ \ \ \ \ \ \ \  $\displaystyle \forall
 R,\ 0<R\leq 1,\ \forall u\in W^{1,r}(B_{R}),\ {\left\Vert{u}\right\Vert}_{L^{t}(B_{R})}\leq
 CR^{-1}\ {\left\Vert{u}\right\Vert}_{W^{1,r}(B_{R})}$\vspace{.5em}  \par 
the constant $C$ depending only on $\displaystyle n,\ r.$
\end{lem}
\quad Proof.\ \par 
Start with $R=1,$ then we have, by the Sobolev embeddings which
 are valid on balls in ${\mathbb{R}}^{n},$ with $t::\frac{1}{t}=\frac{1}{r}-\frac{1}{n},$\
 \par 
\quad \quad \quad \begin{equation} \forall v\in W^{1,r}(B_{1}),\ {\left\Vert{v}\right\Vert}_{L^{t}(B_{1})}\leq
 C{\left\Vert{v}\right\Vert}_{W^{1,r}(B_{1})}\label{aS0}\end{equation}\ \par 
where $\displaystyle C$ depends only on $n$ and $\displaystyle
 r.$ For $u\in W^{1,r}(B_{R})$ we set\ \par 
\vspace{.5em} \ \ \ \ \ \ \ \ \ \ \ \  $\displaystyle \forall
 x\in B_{1},\ y:=Rx\in B_{R},\ v(x):=u(y).$\vspace{.5em}  \ \par 
Then we have\ \par 
\vspace{.5em} \ \ \ \ \ \ \ \ \ \ \ \  $\displaystyle \partial
 v(x)=\partial u(y){\times}\frac{\partial y}{\partial x}=R\partial
 u(y);$\vspace{.5em}  \ \par 
So we get, because the jacobian for this change of variables is $R^{-n},$\ \par 
\vspace{.5em} \ \ \ \ \ \ \ \ \ \ \ \  $\displaystyle {\left\Vert{\partial
 v}\right\Vert}_{L^{r}(B_{1})}^{r}=\int_{B_{1}}{\left\vert{\partial
 v(x)}\right\vert ^{r}dm(x)}=\int_{B_{R}}{\left\vert{\partial
 u(y)}\right\vert ^{r}\frac{R^{r}}{R^{n}}dm(x)}=R^{r-n}{\left\Vert{\partial
 u}\right\Vert}_{L^{r}(B_{R})}^{r}.$\vspace{.5em}  \ \par 
So\ \par 
\quad \quad \quad \begin{equation} \ {\left\Vert{\partial u}\right\Vert}_{L^{r}(B_{R})}=R^{-1+n/r}{\left\Vert{\partial
 v}\right\Vert}_{L^{r}(B_{1})}.\label{HS50}\end{equation}\ \par 
And of course ${\left\Vert{u}\right\Vert}_{L^{r}(B_{R})}=R^{n/r}{\left\Vert{v}\right\Vert}_{L^{r}(B_{1})}.$\
 \par 
So with~\ref{aS0} we get\ \par 
\quad \quad \quad \begin{equation} \ {\left\Vert{u}\right\Vert}_{L^{t}(B_{R})}=R^{n/t}{\left\Vert{v}\right\Vert}_{L^{t}(B_{1})}\leq
 CR^{n/t}{\left\Vert{v}\right\Vert}_{W^{1,r}(B_{1})}.\label{3S1}\end{equation}\
 \par 
But\ \par 
\vspace{.5em} \ \ \ \ \ \ \ \ \ \ \ \ {\it  }$\displaystyle {\left\Vert{u}\right\Vert}_{W^{1,r}(B_{R})}={\left\Vert{u}\right\Vert}_{L^{r}(B_{R})}+{\left\Vert{\partial
 u}\right\Vert}_{L^{r}(B_{R})},$\vspace{.5em} {\it  }\ \par 
and\ \par 
\vspace{.5em} \ \ \ \ \ \ \ \ \ \ \ \  $\displaystyle {\left\Vert{v}\right\Vert}_{W^{1,r}(B_{1})}={\left\Vert{v}\right\Vert}_{L^{r}(B_{1})}+{\left\Vert{\partial
 v}\right\Vert}_{L^{r}(B_{1})},$\vspace{.5em}  \ \par 
so\ \par 
\vspace{.5em} \ \ \ \ \ \ \ \ \ \ \ \  $\displaystyle {\left\Vert{v}\right\Vert}_{W^{1,r}(B_{1})}:=R^{-n/r}{\left\Vert{u}\right\Vert}_{L^{r}(B_{R})}+R^{1-n/r}{\left\Vert{\partial
 u}\right\Vert}_{L^{r}(B_{R})}.$\vspace{.5em}  \ \par 
Because we have $\displaystyle R\leq 1,$ we get\ \par 
\vspace{.5em} \ \ \ \ \ \ \ \ \ \ \ \  $\displaystyle {\left\Vert{v}\right\Vert}_{W^{1,r}(B_{1})}\leq
 R^{-n/r}({\left\Vert{u}\right\Vert}_{L^{r}(B_{R})}+{\left\Vert{\partial
 u}\right\Vert}_{L^{r}(B_{R})})=R^{-n/r}{\left\Vert{u}\right\Vert}_{W^{1,r}(B_{R})}.$\vspace{.5em}
  \ \par 
Putting it in~(\ref{3S1}) we get\ \par 
\vspace{.5em} \ \ \ \ \ \ \ \ \ \ \ \  $\displaystyle {\left\Vert{u}\right\Vert}_{L^{t}(B_{R})}\leq
 CR^{n/t}{\left\Vert{v}\right\Vert}_{W^{1,r}(B_{1})}\leq CR^{-n(\frac{1}{r}-\frac{1}{t})}{\left\Vert{u}\right\Vert}_{W^{1,r}(B_{R})}.$\vspace{.5em}
  \ \par 
But, because $t::\frac{1}{t}=\frac{1}{r}-\frac{1}{n},$ we get
 $(\frac{1}{r}-\frac{1}{t})=\frac{1}{n}$ and\ \par 
\vspace{.5em} \ \ \ \ \ \ \ \ \ \ \ \  $\displaystyle {\left\Vert{u}\right\Vert}_{L^{t}(B_{R})}\leq
 CR^{-1}{\left\Vert{u}\right\Vert}_{W^{1,r}(B_{R})}.$\vspace{.5em}  \ \par 
\quad The constant $C$ depends only on $n,r.$ The proof is complete.
 $\hfill\blacksquare $\ \par 

\begin{lem}
~\label{3S2}Let $x\in M$ and $B(x,R)$ be a $(0,\epsilon )$-admissible
 ball in the complete riemannian manifold $(M,g)$; we have, with
 $1/s=1/r-1/n,$\par 
\vspace{.5em} \ \ \ \ \ \ \ \ \ \ \ \  $\displaystyle \forall
 u\in W^{1,r}(B(x,R)),\ {\left\Vert{u}\right\Vert}_{L^{s}(B(x,R))}\leq
 CR^{-2}{\left\Vert{u}\right\Vert}_{W^{1,r}(B(x,R))},$\vspace{.5em}  \par 
the constant $\displaystyle C$ depending only on $n,\ r$ and
 $\displaystyle \epsilon .$
\end{lem}
\quad Proof.\ \par 
By Lemma~\ref{3S3}, we get in ${\mathbb{R}}^{n}$:\ \par 
\vspace{.5em} \ \ \ \ \ \ \ \ \ \ \ \  $\displaystyle {\left\Vert{u}\right\Vert}_{L^{s}(B_{R})}\leq
 CR^{-1}\ {\left\Vert{u}\right\Vert}_{W^{1,r}(B_{R})}$\vspace{.5em}  \ \par 
so we can apply the comparison Lemma~\ref{HCS43}:\ \par 
\vspace{.5em} \ \ \ \ \ \ \ \ \ \ \ \  $\displaystyle {\left\Vert{u}\right\Vert}_{L^{s}(B(x,R))}\leq
 C{\left\Vert{u\circ \varphi ^{-1}}\right\Vert}_{L^{s}(B_{R}))},$\vspace{.5em}
  \ \par 
which gives:\ \par 
\vspace{.5em} \ \ \ \ \ \ \ \ \ \ \ \  $\displaystyle {\left\Vert{u}\right\Vert}_{L^{s}(B(x,R))}\leq
 CR^{-1}\ {\left\Vert{u}\right\Vert}_{W^{1,r}(B_{R})}.$\vspace{.5em}  \ \par 
Again with the reverse inequalities in the comparison Lemma~\ref{HCS43}:\ \par 
\vspace{.5em} \ \ \ \ \ \ \ \ \ \ \ \  $\displaystyle {\left\Vert{u}\right\Vert}_{W^{1,r}(B_{R})}\leq
 CR^{-1}{\left\Vert{u}\right\Vert}_{W^{1,r}(B(x,R))}.$\vspace{.5em}  \ \par 
So  we get\ \par 
\vspace{.5em} \ \ \ \ \ \ \ \ \ \ \ \  $\displaystyle {\left\Vert{u}\right\Vert}_{L^{s}(B(x,R))}\leq
 CR^{-2}\ {\left\Vert{u}\right\Vert}_{W^{1,r}(B(x,R)}.$\vspace{.5em}  \ \par 
The constant $C$ being independent of $x\in M$ and of $R.$ The
 proof is complete. $\hfill\blacksquare $\ \par 

\subsection{Sobolev comparison estimates for sections of $G.$}

\begin{lem}
~\label{SC19}We have the Sobolev comparison estimates where $B(x,R)$
 is a $(1,\epsilon )$-admissible ball in $M$ and $\varphi \ :\
 B(x,R)\rightarrow {\mathbb{R}}^{n}$  is the admissible chart
 relative to $B(x,R).$ Set $v:=\varphi ^{*}\omega ,$ then:\par 
\vspace{.5em} \ \ \ \ \ \ \ \ \ \ \ \  $\displaystyle \forall
 \omega \in W_{G}^{1,r}(B(x,R)),\ {\left\Vert{\omega }\right\Vert}_{W_{G}^{1,r}(B(x,R))}\leq
 (1+C\epsilon )R^{-1}{\left\Vert{v}\right\Vert}_{W^{1,r}(\varphi
 (B(x,R)))},$\vspace{.5em}  \par 
and, with $B_{e}(0,t)$ the euclidean ball in ${\mathbb{R}}^{n}$
 centered at $0$ and of radius $t,$\par 
\vspace{.5em} \ \ \ \ \ \ \ \ \ \ \ \  $\displaystyle {\left\Vert{v}\right\Vert}_{W^{1,r}(B_{e}(0,(1-\epsilon
 )R))}\leq (1+C\epsilon )R^{-1}{\left\Vert{\omega }\right\Vert}_{W^{1,r}(B(x,R))}.$\vspace{.5em}
  \par 
\quad We also have:\par 
\vspace{.5em} \ \ \ \ \ \ \ \ \ \ \ \  $\displaystyle \forall
 \omega \in L_{G}^{r}(B(x,R)),\ {\left\Vert{\omega }\right\Vert}_{L_{G}^{r}(B(x,R))}\leq
 (1+C\epsilon ){\left\Vert{v}\right\Vert}_{L^{r}(\varphi (B(x,R)))},$\vspace{.5em}
  \par 
and\par 
\vspace{.5em} \ \ \ \ \ \ \ \ \ \ \ \  $\displaystyle {\left\Vert{v}\right\Vert}_{L^{r}(B_{e}(0,(1-\epsilon
 )R))}\leq (1+C\epsilon ){\left\Vert{\omega }\right\Vert}_{L_{G}^{r}(B(x,R))}.$\vspace{.5em}
  
\end{lem}
\quad Proof.\ \par 
We have to compare the norms of $\omega ,\ \nabla \omega ,$ with
 the corresponding ones for $v:=\varphi ^{*}\omega $ in ${\mathbb{R}}^{n}.$\
 \par 
By Lemma~\ref{SB30} the $\epsilon $-admissible ball $B(x,R)$
 trivialises the bundle $G$ hence the image of a section of $G$
 in ${\mathbb{R}}^{n}$ is just vectors of functions. Precisely
 $v:=\varphi ^{*}\omega \in \varphi (B(x,R)){\times}{\mathbb{R}}^{N}.$\ \par 
We have because $(1-\epsilon )\delta _{ij}\leq g_{ij}\leq (1+\epsilon
 )\delta _{ij}$ in $B(x,R)$:\ \par 
\vspace{.5em} \ \ \ \ \ \ \ \ \ \ \ \  $\displaystyle B_{e}(0,(1-\epsilon
 )R)\subset \varphi (B(x,R))\subset B_{e}(0,(1+\epsilon )R).$\vspace{.5em}
  \ \par 
Let $\omega $ be a section of $G$ in $M.$ By our assumption (CMT)
 we have that $\nabla \omega $ depends on the first order derivatives
 of the metric tensor $g.$\ \par 
\quad Because of~(\ref{SB34}) we get, with the fact that $B(x,R)$ is
 $(1,\epsilon )$-admissible, with $\eta :=\frac{\epsilon }{R},$\ \par 
\vspace{.5em} \ \ \ \ \ \ \ \ \ \ \ \  $\displaystyle \sum_{\left\vert{\beta
 }\right\vert =1}{R\sup \ _{i,j=1,...,n,\ y\in B_{x}(R)}\left\vert{\partial
 ^{\beta }g_{ij}(y)}\right\vert }\leq \epsilon \Rightarrow \left\vert{\Gamma
 ^{k}_{ij}}\right\vert \leq C\eta .$\vspace{.5em}  \ \par 
Hence\ \par 
\vspace{.5em} \ \ \ \ \ \ \ \ \ \ \ \  $\displaystyle \forall
 y\in B(x,R),\ \left\vert{\omega (y)}\right\vert =\left\vert{v(z)}\right\vert
 ,\ \ \left\vert{\nabla \omega (y)}\right\vert \leq \left\vert{\partial
 \omega }\right\vert +\left\vert{\Phi }\right\vert ,$\vspace{.5em}  \ \par 
where $\Phi $ depends on the coefficients of $\omega $ and on
 the first order derivatives of the metric tensor $g.$\ \par 
So\ \par 
\quad \quad \quad \begin{equation} \ \left\vert{\nabla \omega (y)}\right\vert \leq
 \left\vert{\partial v(z)}\right\vert +C\eta \left\vert{v(z)}\right\vert
 .\label{SC22}\end{equation}\ \par 
\quad Integrating this we get\ \par 
\vspace{.5em} \ \ \ \ \ \ \ \ \ \ \ \  $\displaystyle {\left\Vert{\nabla
 \omega (y)}\right\Vert}_{L^{r}(B(x,R))}\leq {\left\Vert{\partial
 v}\right\Vert}_{L^{r}(B_{e}(0,(1+\epsilon )R))}+C\eta {\left\Vert{v}\right\Vert}_{L^{r}(B_{e}(0,(1+\epsilon
 )R))}.$\vspace{.5em}  \ \par 
We also have the reverse estimates\ \par 
\vspace{.5em} \ \ \ \ \ \ \ \ \ \ \ \  $\displaystyle {\left\Vert{\partial
 v}\right\Vert}_{L^{r}(B_{e}(0,(1-\epsilon )R))}\leq {\left\Vert{\nabla
 \omega }\right\Vert}_{L^{r}(B(x,R))}+C\eta {\left\Vert{\omega
 }\right\Vert}_{L^{r}(B(x,R))}.$\vspace{.5em}  \ \par 
So, using that\ \par 
\vspace{.5em} \ \ \ \ \ \ \ \ \ \ \ \  $\displaystyle {\left\Vert{\omega
 }\right\Vert}_{W_{G}^{1,r}(B(x,R))}={\left\Vert{\nabla \omega
 }\right\Vert}_{L_{G}^{r}(B(x,R))}+{\left\Vert{\omega }\right\Vert}_{L_{G}^{r}(B(x,R))},$\vspace{.5em}
  \ \par 
we get\ \par 
\vspace{.5em} \ \ \ \ \ \ \ \ \ \ \ \  $\displaystyle {\left\Vert{\omega
 }\right\Vert}_{W_{G}^{1,r}(B(x,R))}\leq {\left\Vert{\partial
 v}\right\Vert}_{L^{r}(B_{e}(0,(1+\epsilon )R))}+C\eta {\left\Vert{v}\right\Vert}_{L^{r}(B_{e}(0,(1+\epsilon
 )R))}\leq $\vspace{.5em}  \ \par 
\vspace{.5em} \ \ \ \ \ \ \ \ \ \ \ \ \ \ \ \ \ \ \ \ \ \ \ \
  $\displaystyle \leq (1+\eta C){\left\Vert{v}\right\Vert}_{W^{1,r}(B_{e}(0,(1+\epsilon
 )R))}.$\vspace{.5em}  \ \par 
Now for $R\leq 1,$ we get:\ \par 
\vspace{.5em} \ \ \ \ \ \ \ \ \ \ \ \  $\displaystyle 1+C\eta
 =1+C\frac{\epsilon }{R}=R^{-1}(R+C\epsilon )\leq R^{-1}(1+C\epsilon
 ).$\vspace{.5em}  \ \par 
The second part follows the same lines but is easier because
 there is no derivations.\ \par 
\quad The proof of the lemma is complete. $\hfill\blacksquare $\ \par 

\begin{lem}
~\label{SC20}Let $x\in M$ and $B(x,R)$ be a $(1,\epsilon )$-admissible
 ball in the complete riemannian manifold $(M,g)$; we have, with
 $1/s=1/r-1/n,$\par 
\vspace{.5em} \ \ \ \ \ \ \ \ \ \ \ \  $\displaystyle \forall
 \omega \in W_{G}^{1,r}(B(x,R)),\ {\left\Vert{\omega }\right\Vert}_{L_{G}^{s}(B(x,R))}\leq
 CR^{-2}{\left\Vert{\omega }\right\Vert}_{W_{G}^{1,r}(B(x,R))}.$\vspace{.5em}
  \par 
the constant $\displaystyle C$ depending only on $n,\ r$ and
 $\displaystyle \epsilon .$
\end{lem}
\quad Proof.\ \par 
Because the image of a section of $G$ in ${\mathbb{R}}^{n}$ is
 just vectors of functions by Lemma~\ref{SB30}, Lemma~\ref{3S3}
 is also true for $\varphi ^{*}\omega $:\ \par 
\vspace{.5em} \ \ \ \ \ \ \ \ \ \ \ \  $\displaystyle {\left\Vert{\varphi
 ^{*}\omega }\right\Vert}_{L^{s}(B_{R})}\leq CR^{-1}\ {\left\Vert{\varphi
 ^{*}\omega }\right\Vert}_{W^{1,r}(B_{R})}$\vspace{.5em}  \ \par 
so we can apply the second part in the comparison Lemma~\ref{SC19}:\ \par 
\vspace{.5em} \ \ \ \ \ \ \ \ \ \ \ \  $\displaystyle {\left\Vert{\omega
 }\right\Vert}_{L_{G}^{s}(B(x,R))}\leq C{\left\Vert{\varphi ^{*}\omega
 }\right\Vert}_{L^{s}(B_{R}))},$\vspace{.5em}  \ \par 
which gives:\ \par 
\vspace{.5em} \ \ \ \ \ \ \ \ \ \ \ \  $\displaystyle {\left\Vert{\omega
 }\right\Vert}_{L_{G}^{s}(B(x,R))}\leq CR^{-1}\ {\left\Vert{\varphi
 ^{*}\omega }\right\Vert}_{W^{1,r}(B_{R})}.$\vspace{.5em}  \ \par 
Again with the reverse inequalities in the comparison Lemma~\ref{SC19}:\ \par 
\vspace{.5em} \ \ \ \ \ \ \ \ \ \ \ \  $\displaystyle {\left\Vert{\varphi
 ^{*}\omega }\right\Vert}_{W^{1,r}(B_{R})}\leq CR^{-1}{\left\Vert{\omega
 }\right\Vert}_{W_{G}^{1,r}(B(x,R))}.$\vspace{.5em}  \ \par 
So  we get\ \par 
\vspace{.5em} \ \ \ \ \ \ \ \ \ \ \ \  $\displaystyle {\left\Vert{\omega
 }\right\Vert}_{L_{G}^{s}(B(x,R))}\leq CR^{-2}\ {\left\Vert{\omega
 }\right\Vert}_{W_{G}^{1,r}(B(x,R)}.$\vspace{.5em}  \ \par 
The constant $C$ is independent of $x\in M$ and of $R.$ The proof
 is complete. $\hfill\blacksquare $\ \par 

\subsection{Local Gaffney type inequality in $L^{r}.$~\label{SB39}}
\quad We shall restrict here to the case of the bundle $\Lambda ^{p}(M)$
 of $p$-forms on $M.$\ \par 
\ \par 
\quad Of course the operator $d$ on $p$-forms is local and so is $\displaystyle
 d^{*}$ as a first order differential operator on $M.$\ \par 
\quad Let $B:=B(x_{0},R)$ be a $(1,\epsilon )$-admissible ball in the
 complete riemannian manifold $(M,g)$ and $(B,\varphi )$ be a
 coordinates chart. Let $\omega $ be a $p$-form in $M.$ Let $\chi
 $ be a smooth cut-off function, $\chi \in {\mathcal{C}}_{0}^{1}(B),\
 0\leq \chi \leq 1,\ \chi \equiv 1$ in $B^{1}:=B(x_{0},R/2).$
 We consider the $p$-form $\chi \omega .$\ \par 
\quad Read in the chart $(B,\varphi )$ with the local coordinates $x,$
 we get $\omega =a_{J}dx^{J}$ with $J=(j_{1},...,j_{p})$ is a
 multi-index of length $p$ and the functions $a_{J}$ are in $W^{1,r}(B).$
 We get then that $d(\chi \omega )=\chi da_{J}\wedge dx^{J}+a_{J}d\chi
 \wedge dx^{J},$ hence, with  $da_{J}=\frac{\partial a_{J}}{\partial
 x^{j}}dx^{j},$ we deduce $d\omega =\frac{\partial a_{J}}{\partial
 x^{j}}dx^{j}\wedge dx^{J}$ and $d(\chi \omega )=\chi \frac{\partial
 a_{J}}{\partial x^{j}}dx^{j}\wedge dx^{J}+a_{J}d\chi \wedge dx^{J}.$\ \par 
\quad We shall take the following notation from the book by C. Voisin~\cite{Voisin02}.\
 \par 
\quad With a local Hodge $*$ operator $\Lambda ^{p}\rightarrow \Lambda
 ^{n-p}$ (locally $M$ is always orientable), i.e. in a coordinates
 chart $U$ with $\omega =a_{J}dx^{J}\in \Lambda ^{p},$ it is defined as:\ \par 
\vspace{.5em} \ \ \ \ \ \ \ \ \ \ \ \  $\displaystyle *\omega
 :=(-1)^{\sigma (\ J)}a_{J^{c}}dx^{J^{c}}\in \Lambda ^{n-p}$\vspace{.5em}
  \ \par 
with $J^{c}$ is the complement of $J$ in $(1,2,...,n)$ and $\sigma
 (J)$ is $0$ or $1.$ We have: $\int_{U}{\omega \wedge \omega
 ^{*}}=\int_{U}{\sum_{J,\ \left\vert{J}\right\vert =p}{}\left\vert{a_{J}}\right\vert
 ^{2}}dv.$\ \par 
Using the link between the $*$ Hodge operator and the adjoint
 $d^{*}$ of $d,$ (see ~\cite[Section 5.1.2, p. 118]{Voisin02}),
 we get: $d^{*}=(-1)^{p}*^{-1}d*$ on $\Lambda ^{p}.$\ \par 
Hence here we have:\ \par 
\vspace{.5em} \ \ \ \ \ \ \ \ \ \ \ \  $\displaystyle d^{*}(\chi
 \omega )=(-1)^{p(n-p)-n+p-1}(*d(*(\chi \omega )))=L_{K}(\frac{\partial
 a_{J}}{\partial x^{j}})dx^{K}+F_{K}(a_{J},\frac{\partial \chi
 }{\partial x^{j}})dx^{K},$\vspace{.5em}  \ \par 
where $L_{K}$ is linear in the $\frac{\partial a_{J}}{\partial
 x^{j}}$ and $F_{K}$ is linear in the $a_{J}$ and in the $\frac{\partial
 \chi }{\partial x^{j}},$ and  $\left\vert{K}\right\vert =p-1.$
 Moreover $L_{K},\ F_{K}$ have compact support in $B.$\ \par 
\ \par 
\quad Now for the covariant derivative $\nabla _{M}$ on $M,$ we have,
 in our chart $(B,\varphi ),$ using~(\ref{SC22})\ \par 
\quad \quad \quad \begin{equation} \left\vert{\nabla _{M}(\chi \omega )}\right\vert
 \leq \left\vert{\partial (\chi \omega )}\right\vert +c\frac{\epsilon
 }{R}\left\vert{\chi \omega }\right\vert .\label{SC23}\end{equation}\ \par 
On the other hand, in ${\mathbb{R}}^{n}$ we have $\left\vert{\nabla
 _{{\mathbb{R}}^{n}}(\chi \omega )}\right\vert =\left\vert{\partial
 (\chi \omega )}\right\vert $ just because $\nabla _{{\mathbb{R}}^{n}}(\chi
 \omega )=\partial (\chi \omega ).$\ \par 
Now we are in position to apply the following Proposition 4.3
 in~\cite{Scott95}:\ \par 

\begin{prop}
~\label{SC24}(Gaffney type inequality for $L^{r}$)\par 
\vspace{.5em} \ \ \ \ \ \ \ \ \ \ \ \  $\displaystyle \ {\left\Vert{\nabla
 _{{\mathbb{R}}^{n}}\omega }\right\Vert}_{r}^{r}\leq C{\left\Vert{d\omega
 }\right\Vert}_{L^{r}({\mathbb{R}}^{n})}^{r}+{\left\Vert{d^{*}\omega
 }\right\Vert}_{L^{r}({\mathbb{R}}^{n})}^{r}$\vspace{.5em}  \par 
\ for $\omega \in {\mathcal{C}}^{\infty }_{0}(\Lambda ^{p}({\mathbb{R}}^{n}))$
 and $C=C(n,r).$
\end{prop}
\quad From it, we get in ${\mathbb{R}}^{n},$\ \par 
\vspace{.5em} \ \ \ \ \ \ \ \ \ \ \ \  $\displaystyle {\left\Vert{\partial
 (\chi \omega )}\right\Vert}_{L^{r}(\varphi (B))}\leq C({\left\Vert{d(\chi
 \omega )}\right\Vert}_{L^{r}(\varphi (B))}+{\left\Vert{d^{*}(\chi
 \omega )}\right\Vert}_{L^{r}(\varphi (B))}).$\vspace{.5em}  \ \par 
So, by~(\ref{SC23}), we get, because $\chi \equiv 1$ in $B^{1}$:\ \par 
\vspace{.5em} \ \ \ \ \ \ \ \ \ \ \ \  $\displaystyle \ {\left\Vert{\nabla
 _{M}\omega }\right\Vert}_{L^{r}(B^{1})}\leq C({\left\Vert{d(\chi
 \omega )}\right\Vert}_{L^{r}(B)}+{\left\Vert{d^{*}(\chi \omega
 )}\right\Vert}_{L^{r}(B)})+c\epsilon R^{-1}{\left\Vert{\chi
 \omega }\right\Vert}_{L^{r}(B)},$\vspace{.5em}  \ \par 
because, by condition (*) in the definition of the $\epsilon
 $-admissible ball $B,$ we have that the Lebesgue measure in
 $\varphi (B)$ and the volume measure on $B$ are equivalent.\ \par 
\quad Now, because ${\left\Vert{d\chi }\right\Vert}_{\infty }\leq cR^{-1}$
 in $B$:\ \par 
\vspace{.5em} \ \ \ \ \ \ \ \ \ \ \ \  $\displaystyle {\left\Vert{d(\chi
 \omega )}\right\Vert}_{L^{r}(B)}\leq {\left\Vert{d\chi }\right\Vert}_{\infty
 }{\left\Vert{\omega }\right\Vert}_{L^{r}(B)}+{\left\Vert{d\omega
 }\right\Vert}_{L^{r}(B)}\leq cR^{-1}{\left\Vert{\omega }\right\Vert}_{L^{r}(B)}+{\left\Vert{d\omega
 }\right\Vert}_{L^{r}(B)},$\vspace{.5em}  \ \par 
and the same\ \par 
\vspace{.5em} \ \ \ \ \ \ \ \ \ \ \ \  $\displaystyle {\left\Vert{d^{*}(\chi
 \omega )}\right\Vert}_{L^{r}(B)}\leq cR^{-1}{\left\Vert{\omega
 }\right\Vert}_{L^{r}(B)}+{\left\Vert{d^{*}\omega }\right\Vert}_{L^{r}(B)}.$\vspace{.5em}
  \ \par 
Hence, with new constants $c,C$ depending only on $n,\ r,\ p,\
 \epsilon ,$\ \par 
\vspace{.5em} \ \ \ \ \ \ \ \ \ \ \ \  $\displaystyle {\left\Vert{\nabla
 _{M}\omega }\right\Vert}_{L^{r}(B^{1})}\leq C({\left\Vert{d(\omega
 )}\right\Vert}_{L^{r}(B)}+{\left\Vert{d^{*}(\omega )}\right\Vert}_{L^{r}(B)})+cR^{-1}{\left\Vert{\omega
 }\right\Vert}_{L^{r}(B)}.$\vspace{.5em}  \ \par 
\quad So we proved the corollary of Scott's Proposition~\ref{SC24}:\ \par 

\begin{cor}
~\label{SC25}Let $B:=B(x_{0},R)$ be a $(1,\epsilon )$-admissible
 ball in the complete riemannian manifold $(M,g)$ and set $B^{1}:=B(x_{0},R/2).$
 Let $\omega $ be a $p$-form in $M.$ We have the local $L^{r}$
 Gaffney's inequality:\par 
\vspace{.5em} \ \ \ \ \ \ \ \ \ \ \ \  $\displaystyle {\left\Vert{\nabla
 _{M}\omega }\right\Vert}_{L^{r}(B^{1})}\leq C({\left\Vert{d(\omega
 )}\right\Vert}_{L^{r}(B)}+{\left\Vert{d^{*}(\omega )}\right\Vert}_{L^{r}(B)})+cR^{-1}{\left\Vert{\omega
 }\right\Vert}_{L^{r}(B)},$\vspace{.5em}  \par 
the constants $c,C$ depending only on $n,\ r,\ p,\ \epsilon .$
\end{cor}

\section{Vitali covering.~\label{SB36}}

\begin{lem}
Let ${\mathcal{F}}$ be a collection of balls $\lbrace B(x,r(x))\rbrace
 $ in a metric space, with $\forall B(x,r(x))\in {\mathcal{F}},\
 0<r(x)\leq R.$ There exists a disjoint subcollection ${\mathcal{G}}$
 of ${\mathcal{F}}$ with the following properties:\par 
\quad \quad every ball $B$ in ${\mathcal{F}}$ intersects a ball $C$ in ${\mathcal{G}}$
 and $B\subset 5C.$
\end{lem}
This is a well known lemma, see for instance ~\cite{EvGar92},
 section 1.5.1.\ \par 
\ \par 
\quad Fix $\epsilon >0$ and let $\forall x\in M,\ r(x):=R_{\epsilon
 }(x)/10,\ $where $R_{\epsilon }(x)$ is the admissible radius
 at $\displaystyle x,$ we built a Vitali covering with the collection
 ${\mathcal{F}}:=\lbrace B(x,r(x))\rbrace _{x\in M}.$ The previous
 lemma gives a disjoint subcollection ${\mathcal{G}}$ such that
 every ball $B$ in ${\mathcal{F}}$ intersects a ball $C$ in ${\mathcal{G}}$
 and we have $\displaystyle B\subset 5C.$ We set ${\mathcal{D}}(\epsilon
 ):=\lbrace x\in M::B(x,r(x))\in {\mathcal{G}}\rbrace $ and ${\mathcal{C}}_{\epsilon
 }:=\lbrace B(x,5r(x)),\ x\in {\mathcal{D}}(\epsilon )\rbrace
 $: we shall call ${\mathcal{C}}_{\epsilon }$ a $\epsilon $-\textbf{admissible
 covering} of $(M,g).$ We notice that $B(x,5r(x))=B(x,R_{\epsilon
 }(x)/2).$\ \par 
\ \par 
\quad Then we have~\cite[Proposition 6.12]{SobPar18},\ \par 

\begin{prop}
Let $(M,g)$ be a riemannian manifold, then the overlap of a $\epsilon
 $-admissible covering ${\mathcal{C}}(\epsilon )$ is less than
 $T=\frac{(1+\epsilon )^{n/2}}{(1-\epsilon )^{n/2}}(100)^{n},$ i.e.\par 
\quad \quad \quad $\forall x\in M,\ x\in B(y,5r(y)),$ where $B(y,r(y))\in {\mathcal{G}}(\epsilon
 ),$ for at most $T$ such balls.\par 
Moreover we have\par 
\vspace{.5em} \ \ \ \ \ \ \ \ \ \ \ \  $\displaystyle \forall
 f\in L^{1}(M),\ \sum_{j\in {\mathbb{N}}}{\int_{B(x_{j},r(x_{j}))}{\left\vert{f(x)}\right\vert
 dv_{g}(x)}}\leq T{\left\Vert{f}\right\Vert}_{L^{1}(M)}.$\vspace{.5em}  
\end{prop}
\quad And its Corollary~\cite[Corollary 6.13]{SobPar18},\ \par 

\begin{cor}
~\label{SC26}Let $(M,g)$ be a complete riemannian manifold. Consider
 the covering by the balls $\lbrace B(x,R_{\epsilon }(x)),\ x\in
 {\mathcal{D}}(\epsilon )\rbrace .$ Then the overlap of the associated
 covering verifies:\par 
\vspace{.5em} \ \ \ \ \ \ \ \ \ \ \ \  $\displaystyle T_{1}\leq
 \frac{(1+\epsilon )^{n/2}}{(1-\epsilon )^{n/2}}(100)^{n}{\times}2^{n}.$\vspace{.5em}
  
\end{cor}

\section{Global results.~\label{SB40}}

\subsection{Global estimates for sections of $G.$}

\begin{lem}
~\label{S5}We have, for any section $f$ of $G$ with $w(x):=R(x)^{\mu
 }$ and $B(x):=B(x,R(x)),$ with $x\in {\mathcal{D}}(\epsilon
 )$ and $R(x):=R_{\epsilon }(x),$ that:\par 
\vspace{.5em} \ \ \ \ \ \ \ \ \ \ \ \  $\displaystyle \forall
 \tau \geq 1,\ {\left\Vert{f}\right\Vert}_{L^{\tau }_{G}(M,\
 w)}^{\tau }\simeq \sum_{x\in {\mathcal{D}}(\epsilon )}{R(x)^{\mu
 }{\left\Vert{f}\right\Vert}_{L^{\tau }_{G}(B(x))}^{\tau }}.$\vspace{.5em}  
\end{lem}
\quad Proof.\ \par 
We consider the function $\left\vert{f}\right\vert .$ We set
 $w(x):=R(x)^{\mu }$ to get an adapted weight. Let $x\in {\mathcal{D}}(\epsilon
 ),$ then $B(x):=B(x,R(x))\in {\mathcal{C}}_{\epsilon }.$\ \par 
\quad We have, because ${\mathcal{C}}_{\epsilon }$ is a covering of $M,$\ \par 
\vspace{.5em} \ \ \ \ \ \ \ \ \ \ \ \  $\displaystyle {\left\Vert{f}\right\Vert}_{L^{\tau
 }(M,w)}^{\tau }:=\int_{M}{\left\vert{f(x)}\right\vert ^{\tau
 }w(x)dv(x)}\leq \sum_{x\in {\mathcal{D}}(\epsilon )}{\int_{B(x)}{\left\vert{f(y)}\right\vert
 ^{\tau }R(y)^{\mu }}dv(y)}.$\vspace{.5em}  \ \par 
Because we have, by Lemma~\ref{S2}, $\forall y\in B,\ R(y)\leq
 2R(x),$ we get\ \par 
\vspace{.5em} \ \ \ \ \ \ \ \ \ \ \ \  $\displaystyle \sum_{x\in
 {\mathcal{D}}(\epsilon )}{\int_{B(x)}{\left\vert{f(y)}\right\vert
 ^{\tau }R(y)^{\mu }}dv(y)}\leq $\vspace{.5em}  \ \par 
\vspace{.5em} \ \ \ \ \ \ \ \ \ \ \ \ \ \ \ \ \ \ \ \ \ \ \ \
  $\displaystyle \leq \sum_{x\in {\mathcal{D}}(\epsilon )}{2^{\mu
 }R(x)^{\mu }\int_{B(x)}{\left\vert{f(y)}\right\vert ^{\tau }}dv(y)}\leq
 2^{\mu }\sum_{x\in {\mathcal{D}}(\epsilon )}{R(x)^{\mu }{\left\Vert{f}\right\Vert}_{L^{\tau
 }(B(x))}^{\tau }}.$\vspace{.5em}  \ \par 
Hence we get\ \par 
\vspace{.5em} \ \ \ \ \ \ \ \ \ \ \ \  $\displaystyle {\left\Vert{f}\right\Vert}_{L^{\tau
 }(M,w)}^{\tau }\leq 2^{\mu }\sum_{x\in {\mathcal{D}}(\epsilon
 )}{R(x)^{\mu }{\left\Vert{f}\right\Vert}_{L^{\tau }(B)}^{\tau
 }}.$\vspace{.5em}  \ \par 
\ \par 
\quad To get the reverse inequality we still use Lemma~\ref{S2}, to
 have $\forall y\in B,\ R(x)\leq 2R(y)$ so we get:\ \par 
\vspace{.5em} \ \ \ \ \ \ \ \ \ \ \ \  $\displaystyle \sum_{x\in
 {\mathcal{D}}(\epsilon )}{R(x)^{\mu }\int_{B(x)}{\left\vert{f(y)}\right\vert
 ^{\tau }}dv(y)}\leq 2^{\mu }\sum_{x\in {\mathcal{D}}(\epsilon
 )}{\int_{B(x)}{R(y)^{\mu }\left\vert{f(y)}\right\vert ^{\tau
 }}dv(y)}.$\vspace{.5em}  \ \par 
Now we use the fact that the overlap of ${\mathcal{C}}_{\epsilon
 }$ is bounded by $T,$\ \par 
\vspace{.5em} \ \ \ \ \ \ \ \ \ \ \ \  $\displaystyle \sum_{x\in
 {\mathcal{D}}(\epsilon )}{\int_{B(x)}{R(y)^{\mu }\left\vert{f(y)}\right\vert
 ^{\tau }}dv(y)}\leq 2^{\mu }T\int_{M}{R(y)^{\mu }\left\vert{f(y)}\right\vert
 ^{\tau }}dv(y)=2^{\mu }T{\left\Vert{f}\right\Vert}_{L^{\tau
 }(M,w)}^{\tau }.$\vspace{.5em}  \ \par 
So we get\ \par 
\vspace{.5em} \ \ \ \ \ \ \ \ \ \ \ \  $\displaystyle \sum_{x\in
 {\mathcal{D}}(\epsilon )}{R(x)^{\mu }{\left\Vert{f}\right\Vert}_{L^{\tau
 }(B)}}^{\tau }\leq 2^{\mu }T{\left\Vert{f}\right\Vert}_{L^{\tau
 }(M,w)}^{\tau }.$\vspace{.5em}  \ \par 
The proof is complete. $\hfill\blacksquare $\ \par 

\begin{lem}
~\label{S6}Let $u\in W^{1,r}_{G}(M,w)$ where $w(x)=R(x)^{\gamma
 }.$ We have:\par 
\vspace{.5em} \ \ \ \ \ \ \ \ \ \ \ \  $\displaystyle \sum_{B\in
 {\mathcal{C}}_{\epsilon },\ j=0,1}{R(x)^{\gamma }{\left\Vert{\nabla
 ^{j}u}\right\Vert}^{r}_{L^{r}_{G}(B)}}\simeq {\left\Vert{u}\right\Vert}^{r}_{W^{1,r}_{G}(M,w)}$\vspace{.5em}
  \par 
and\par 
\vspace{.5em} \ \ \ \ \ \ \ \ \ \ \ \  $\displaystyle \sum_{B\in
 {\mathcal{C}}_{\epsilon }}{R(x)^{\gamma }{\left\Vert{u}\right\Vert}^{r}_{W^{1,r}_{G}(B)}}\simeq
 {\left\Vert{u}\right\Vert}^{r}_{W^{1,r}_{G}(M,w)}.$\vspace{.5em}  
\end{lem}
\quad Proof.\ \par 
We start with $u\in W^{1,r}_{G}(M,w)$ where $w(x)=R(x)^{\gamma
 }.$ This means:\ \par 
\vspace{.5em} \ \ \ \ \ \ \ \ \ \ \ \  $\displaystyle \forall
 j\in {\mathbb{N}},\ j\leq 1,\ \int_{M}{\left\vert{\nabla ^{j}u(y)}\right\vert
 ^{r}R(y)^{\gamma }dv(y)}\leq C{\left\Vert{u}\right\Vert}^{r}_{W^{1,r}_{p}(M,w)}.$\vspace{.5em}
  \ \par 
\quad Take a Vitali $\epsilon $-covering ${\mathcal{C}}_{\epsilon }$
 of $M.$ Let $B:=B(x,R(x))\in {\mathcal{C}}_{\epsilon }$; because
 $u\in W^{1,r}_{G}(M,w),$ we have $u\in W_{G}^{1,r}(B,w).$\ \par 
By Lemma~\ref{S5} we have, with $f:=\left\vert{\nabla ^{j}u}\right\vert
 $ and $\tau =r,$\ \par 
\vspace{.5em} \ \ \ \ \ \ \ \ \ \ \ \  $\displaystyle \int_{M}{\left\vert{\nabla
 ^{j}u(y)}\right\vert ^{r}R(y)^{\gamma }dv(y)}\simeq \sum_{x\in
 {\mathcal{D}}(\epsilon )}{R(x)^{\gamma }{\left\Vert{\nabla ^{j}u}\right\Vert}_{L^{r}_{G}(B(x))}^{r}}.$\vspace{.5em}
  \ \par 
Hence, adding on $j,\ j\leq 1,$\ \par 
\quad \quad \begin{equation} \sum_{x\in {\mathcal{D}}(\epsilon ),\ j\leq
 1,}{R(x)^{\gamma }{\left\Vert{\nabla ^{j}u}\right\Vert}_{L^{r}_{G}(B(x))}^{r}}\leq
 C{\left\Vert{u}\right\Vert}^{r}_{W^{1,r}_{G}(M,w)}.\label{S3}\end{equation}\
 \par 
Using that ${\left\Vert{u}\right\Vert}_{W^{1,r}(B(x,R))}={\left\Vert{\nabla
 u}\right\Vert}_{L^{r}(B(x,R))}+{\left\Vert{u}\right\Vert}_{L^{r}(B(x,R))},$
 we deduce from~(\ref{S3}):\ \par 
\vspace{.5em} \ \ \ \ \ \ \ \ \ \ \ \  $\displaystyle \sum_{B\in
 {\mathcal{C}}_{\epsilon }}{R(x)^{\gamma }{\left\Vert{u}\right\Vert}^{r}_{W^{1,r}_{G}(B)}}\leq
 C{\left\Vert{u}\right\Vert}^{r}_{W^{1,r}_{G}(M,w)}.$\vspace{.5em}  \ \par 
\ \par 
\quad For the converse, because:\ \par 
\vspace{.5em} \ \ \ \ \ \ \ \ \ \ \ \  $\displaystyle {\left\Vert{u}\right\Vert}^{r}_{W^{1,r}_{G}(M,w)}\leq
 \sum_{j=0}^{1}{\int_{M}{\left\vert{\nabla ^{j}u(y)}\right\vert
 ^{r}R(y)^{\gamma }dv(y)}}$\vspace{.5em}  \ \par 
we get again with Lemma~\ref{S5},\ \par 
\vspace{.5em} \ \ \ \ \ \ \ \ \ \ \ \  $\displaystyle \int_{M}{\left\vert{\nabla
 ^{j}u(y)}\right\vert ^{r}R(y)^{\gamma }dv(y)}\simeq \sum_{x\in
 {\mathcal{D}}(\epsilon )}{R(x)^{\gamma }{\left\Vert{\nabla ^{j}u}\right\Vert}_{L^{r}_{G}(B(x))}^{r}}.$\vspace{.5em}
  \ \par 
Hence\ \par 
\vspace{.5em} \ \ \ \ \ \ \ \ \ \ \ \  $\displaystyle {\left\Vert{u}\right\Vert}^{r}_{W^{1,r}_{G}(M,w)}\leq
 \sum_{x\in {\mathcal{D}}(\epsilon ),\ j=0,1}{R(x)^{\gamma }{\left\Vert{\nabla
 ^{j}u}\right\Vert}_{L^{r}_{G}(B(x))}^{r}}\lesssim \sum_{x\in
 {\mathcal{D}}(\epsilon )}{R(x)^{\gamma }{\left\Vert{u}\right\Vert}_{W^{1,r}_{G}(B(x))}^{r}}.$\vspace{.5em}
  \ \par 
Now from\ \par 
\quad \quad \quad \begin{equation} \sum_{B\in {\mathcal{C}}_{\epsilon },\ j=0,1}{R(x)^{\gamma
 }{\left\Vert{\nabla ^{j}u}\right\Vert}^{r}_{L^{r}_{G}(B)}}\simeq
 {\left\Vert{u}\right\Vert}^{r}_{W^{1,r}_{G}(M,w)}\label{S9}\end{equation}\
 \par 
the constants underlying $\simeq $ being independent of $B,$ and with:\ \par 
\vspace{.5em} \ \ \ \ \ \ \ \ \ \ \ \  $\displaystyle \sum_{\
 j=0,1}{{\left\Vert{\nabla ^{j}u}\right\Vert}^{r}_{L^{r}_{G}(B)}}={\left\Vert{u}\right\Vert}^{r}_{W^{1,r}_{G}(B)},$\vspace{.5em}
  \ \par 
we get\ \par 
\vspace{.5em} \ \ \ \ \ \ \ \ \ \ \ \  $\displaystyle \sum_{\
 j=0,1}{R(x)^{\gamma }{\left\Vert{\nabla ^{j}u}\right\Vert}^{r}_{L^{r}_{G}(B)}}\simeq
 R(x)^{\gamma }{\left\Vert{u}\right\Vert}^{r}_{W^{1,r}_{G}(B)}.$\vspace{.5em}
  \ \par 
So, porting in~(\ref{S9}), we get\ \par 
\vspace{.5em} \ \ \ \ \ \ \ \ \ \ \ \  $\displaystyle \sum_{B\in
 {\mathcal{C}}_{\epsilon },}{R(x)^{\gamma }{\left\Vert{u}\right\Vert}^{r}_{W^{1,r}_{G}(B)}}\simeq
 {\left\Vert{u}\right\Vert}^{r}_{W^{1,r}_{G}(M,w)}.$\vspace{.5em}  \ \par 
The proof is complete. $\hfill\blacksquare $ \ \par 

\begin{lem}
~\label{S10}We have:\par 
\vspace{.5em} \ \ \ \ \ \ \ \ \ \ \ \  $\displaystyle r\geq s\Rightarrow
 \sum_{j\in {\mathbb{N}}}{a_{j}^{r}}\leq {\left({\sum_{j\in {\mathbb{N}}}{a_{j}^{s}}}\right)}^{r/s};\
 \ r\leq s\Rightarrow \sum_{j\in {\mathbb{N}}}{a_{j}^{r}}\geq
 {\left({\sum_{j\in {\mathbb{N}}}{a_{j}^{s}}}\right)}^{r/s}.$\vspace{.5em}  
\end{lem}
\quad Proof.\ \par 
The comparison of the norms of $\ell ^{r}({\mathbb{N}})$ and
 $\displaystyle \ell ^{s}({\mathbb{N}})$ gives the result. $\hfill\blacksquare
 $\ \par 

\begin{lem}
~\label{S8}Let $B=B(x,R)$ be an $\epsilon $-admissible ball in
 $M.$ We have, for $\omega \in L^{s}_{G}(B)$ with $s\geq r$:\par 
\vspace{.5em} \ \ \ \ \ \ \ \ \ \ \ \  $\displaystyle {\left\Vert{\omega
 }\right\Vert}_{L_{G}^{r}(B))}\leq c(n,\epsilon )R^{\frac{n}{r}-\frac{n}{s}}{\left\Vert{\omega
 }\right\Vert}_{L_{G}^{s}(B)},$\vspace{.5em}  \par 
with $c$ depending only on $n,\ \epsilon $ and $G.$
\end{lem}
\quad Proof.\ \par 
First suppose that $B$ is in ${\mathbb{R}}^{n}$ with the Lebesgue
 measure. Let $\omega \in L^{s}(B)).$ Because $\frac{dv}{\left\vert{B}\right\vert
 }$ is a probability measure on $B,$ where $\left\vert{B}\right\vert
 $ is the volume of the ball $B$ and $s\geq r,$ we get\ \par 
\vspace{.5em} \ \ \ \ \ \ \ \ \ \ \ \  $\displaystyle {\left({\int_{B}{\left\vert{\omega
 (y)}\right\vert ^{r}\frac{dv(y)}{\left\vert{B}\right\vert }}}\right)}^{1/r}\leq
 {\left({\int_{B}{\left\vert{\omega (y)}\right\vert ^{s}\frac{dv(y)}{\left\vert{B}\right\vert
 }}}\right)}^{1/s},$\vspace{.5em}  \ \par 
hence\ \par 
\vspace{.5em} \ \ \ \ \ \ \ \ \ \ \ \  $\displaystyle {\left\Vert{\omega
 }\right\Vert}_{L^{r}(B)}\leq \left\vert{B}\right\vert ^{\frac{1}{r}-\frac{1}{s}}{\left\Vert{\omega
 }\right\Vert}_{L^{s}(B)}.$\vspace{.5em}  \ \par 
Now back to the manifold, for $B_{x}:=B(x,R)$ a $\epsilon $-admissible
 ball we get\ \par 
\vspace{.5em} \ \ \ \ \ \ \ \ \ \ \ \  $\displaystyle \forall
 y\in B_{x},\ (1-\epsilon )^{n}\leq \left\vert{\mathrm{d}\mathrm{e}\mathrm{t}g(y)}\right\vert
 \leq (1+\epsilon )^{n},$\vspace{.5em}  \ \par 
hence we have, comparing the Lebesgue measure in ${\mathbb{R}}^{n}$
 with the volume measure in $M,$\ \par 
\vspace{.5em} \ \ \ \ \ \ \ \ \ \ \ \  $\displaystyle \forall
 x\in M,\ (1-\epsilon )^{n/2}\nu _{n}R^{n}\leq \mathrm{V}\mathrm{o}\mathrm{l}(B(x,\
 R_{\epsilon }(x)))\leq (1+\epsilon )^{n/2}\nu _{n}R^{n},$\vspace{.5em}  \ \par 
where $\nu _{n}$ is the volume of the unit ball in ${\mathbb{R}}^{n}.$
 So, because the ball $B$ trivialises the bundle $G$ by Lemma~\ref{SB30},
 on the manifold $M,$ we have\ \par 
\vspace{.5em} \ \ \ \ \ \ \ \ \ \ \ \  $\displaystyle {\left\Vert{\omega
 }\right\Vert}_{L_{G}^{r}(B))}\leq c(n,\epsilon ,G)R^{\frac{n}{r}-\frac{n}{s}}{\left\Vert{\omega
 }\right\Vert}_{L_{G}^{s}(B)}$\vspace{.5em}  \ \par 
with $c$ depending only on $n,\ \epsilon $ and $G.$ $\hfill\blacksquare $\ \par 
\ \par 
\quad Now we fix $x\in {\mathcal{D}}(\epsilon )$ hence the ball $B:=B(x,R_{\epsilon
 }(x))$ is fixed. We have, with $1/s=1/r-1/n,$ that $W^{1,r}(B)\subset
 L^{s}(B),$  by Lemma~\ref{3S2} for functions or by Lemma~\ref{SC20}
 for sections of $G$ with:\ \par 
\vspace{.5em} \ \ \ \ \ \ \ \ \ \ \ \  $\displaystyle \forall
 u\in W_{G}^{1,r}(B(x,R)),\ {\left\Vert{u}\right\Vert}_{L_{G}^{s}(B(x,R))}\leq
 CR(x)^{-2}{\left\Vert{u}\right\Vert}_{W_{G}^{1,r}(B(x,R))}.$\vspace{.5em}
  \ \par 
So\ \par 
\quad \quad \quad \begin{equation} R(x)^{2+\gamma /r}{\left\Vert{u}\right\Vert}_{L_{G}^{s}(B(x,R))}\leq
 CR(x)^{\gamma /r}{\left\Vert{u}\right\Vert}_{W_{G}^{1,r}(B(x,R))}.\label{S7}\end{equation}\
 \par 
\quad Now the aim is to get the global form of the previous result.\ \par 
\ \par 
\quad Let $w'(x):=R(x)^{\nu }$ with $\nu :=s(2+\gamma /r).$ Using~(\ref{S7}),
 raising to the power $s$ and adding, we get:\ \par 
\quad \quad \quad \begin{equation} \sum_{B\in {\mathcal{C}}_{\epsilon }}{R(x)^{\nu
 }{\left\Vert{u}\right\Vert}^{s}_{L^{s}_{G}(B)}}\leq C\sum_{B\in
 {\mathcal{C}}_{\epsilon }}{R(x)^{s\gamma /r}{\left\Vert{u}\right\Vert}^{s}_{W_{G}^{1,r}(B(x,R))}}.\label{SC21}\end{equation}\
 \par 
Because $s\geq r,$ the Lemma~\ref{S10} gives:\ \par 
\vspace{.5em} \ \ \ \ \ \ \ \ \ \ \ \  $\displaystyle \sum_{B\in
 {\mathcal{C}}_{\epsilon }}{R(x)^{s\gamma /r}{\left\Vert{u}\right\Vert}^{s}_{W^{1,r}(B(x,R))}}\leq
 {\left({\sum_{B\in {\mathcal{C}}_{\epsilon }}{R(x)^{\gamma }{\left\Vert{u}\right\Vert}^{r}_{W^{1,r}(B(x,R))}}}\right)}^{s/r}.$\vspace{.5em}
  \ \par 
Now we use Lemma~\ref{S6} with $w(x):=R(x)^{\gamma },$\ \par 
\vspace{.5em} \ \ \ \ \ \ \ \ \ \ \ \  $\displaystyle \sum_{B\in
 {\mathcal{C}}_{\epsilon }}{R(x)^{\gamma }{\left\Vert{u}\right\Vert}^{r}_{W^{1,r}_{G}(B)}}\simeq
 {\left\Vert{u}\right\Vert}^{r}_{W^{1,r}_{G}(M,w)}$\vspace{.5em}  \ \par 
to get, with Lemma~\ref{S5}, with $\mu =\nu :=2s+s\gamma /r,$\ \par 
\vspace{.5em} \ \ \ \ \ \ \ \ \ \ \ \  $\displaystyle {\left\Vert{u}\right\Vert}^{s}_{L^{s}_{G}(M,w')}\simeq
 \sum_{B\in {\mathcal{C}}_{\epsilon }}{R(x)^{\nu }{\left\Vert{u}\right\Vert}^{s}_{L_{G}^{s}(B(x,R))}}$\vspace{.5em}
  \ \par 
hence, by~(\ref{SC21}),\ \par 
\vspace{.5em} \ \ \ \ \ \ \ \ \ \ \ \  $\displaystyle {\left\Vert{u}\right\Vert}^{s}_{L^{s}_{G}(M,w')}\leq
 C{\left({\sum_{B\in {\mathcal{C}}_{\epsilon }}{R(x)^{\gamma
 }{\left\Vert{u}\right\Vert}^{r}_{W_{G}^{1,r}(B(x,R))}}}\right)}^{s/r}\simeq
 C{\left\Vert{u}\right\Vert}^{s}_{W^{1,r}_{G}(M,w)}.$\vspace{.5em}  \ \par 
\quad So we proved the following Sobolev embedding Theorem with weights:\ \par 

\begin{thm}
~\label{SC15}Let $1/s=1/r-1/n>0.$ Let $u\in W^{1,r}_{G}(M,w)$
 with $w(x):=R(x)^{\gamma }.$ Then we have $u\in L^{s}_{G}(M,w'),$
 with $\nu :=s(2+\gamma /r)$ and $w':=R(x)^{\nu },$ with the control:\par 
\vspace{.5em} \ \ \ \ \ \ \ \ \ \ \ \  $\displaystyle {\left\Vert{u}\right\Vert}_{L^{s}_{G}(M,w')}\leq
 C{\left\Vert{u}\right\Vert}_{W^{1,r}_{G}(M,w)}.$\vspace{.5em}  
\end{thm}
\quad Theorem~\ref{SC15} with Proposition~\ref{SC14} give:\ \par 

\begin{thm}
~\label{SB32}Let $(M,g)$ be a complete riemannian manifold. Let
 $G:=(H,\pi ,M)$ be a complex smooth adapted vector bundle over
 $M.$ Let $w(x):=R(x)^{\gamma }$ and $w':=R(x)^{\nu }$ with 
 $\nu :=s(2+\gamma /r).$ Then $W_{G}^{m,r}(M,w)$ is embedded
 in $W_{G}^{k,s}(M,w'),$ with $\frac{1}{s}=\frac{1}{r}-\frac{(m-k)}{n}>0$
 and:\par 
\vspace{.5em} \ \ \ \ \ \ \ \ \ \ \ \  $\displaystyle \forall
 u\in W_{G}^{m,r}(M,w),\ {\left\Vert{u}\right\Vert}_{W_{G}^{k,s}(M,w')}\leq
 C{\left\Vert{u}\right\Vert}_{W_{G}^{m,r}(M,w)}.$\vspace{.5em}  
\end{thm}

\begin{rem}
The weights, as function of the $\epsilon $-admissible radius
 $R_{\epsilon }(x),$ do not depend on the fact that we work with
 functions or sections of $G$ but  the radius itself depends
 on that fact. The $(1,\epsilon )$ admissible radius for the
 sections of $G$ is smaller than the $(0,\epsilon )$ one for functions.
\end{rem}
\ \par 
\quad In ~\cite[Theorem 6.23, p. 21]{SobPar18}, we proved, with the
 bundle $\Lambda ^{p}(M)$ of $p$-forms, the following:\ \par 

\begin{thm}
~\label{SB31}Let $M$ be a complete {\bf non compact} riemannian
 manifold of class ${\mathcal{C}}^{2}$ without boundary. Let
 $\alpha >0,\ r\geq 2$ and $k$ the smallest integer such that,
 with $\frac{1}{r_{k}}=\frac{1}{2}-\frac{2k}{n},$ we have $r_{k}\geq
 r.$ For any $\omega \in L^{r}(\lbrack 0,T+\alpha \rbrack ,L^{r}_{p}(M,w_{1}))\cap
 L^{r}(\lbrack 0,T+\alpha \rbrack ,L^{2}_{p}(M))$ there is a
 $u\in L^{r}(\lbrack 0,T\rbrack ,W^{2,r}_{p}(M,w_{2}))$ such
 that $\partial _{t}u+\Delta u=\omega $ and:\par 
\vspace{.5em} \ \ \ \ \ \ \ \ \ \ \ \  $\displaystyle {\left\Vert{\partial
 _{t}u}\right\Vert}_{L^{r}(\lbrack 0,T\rbrack ,L^{r}_{p}(M,w_{2}))}+{\left\Vert{u}\right\Vert}_{L^{r}(\lbrack
 0,T\rbrack ,W^{2,r}_{p}(M,w_{2}))}\leq $\vspace{.5em}  \par 
\vspace{.5em} \ \ \ \ \ \ \ \ \ \ \ \ \ \ \ \ \ \ \ \ \ \ \ \
  $\displaystyle \leq c_{1}{\left\Vert{\omega }\right\Vert}_{L^{r}(\lbrack
 0,T+\alpha \rbrack ,L^{r}_{p}(M,w_{1}))}+c_{2}{\sqrt{T+\alpha
 }}{\left\Vert{\omega }\right\Vert}_{L^{r}(\lbrack 0,T+\alpha
 \rbrack ,L^{2}_{p}(M))}$\vspace{.5em}  \par 
where the weights functions are: $w_{1}(x)=R(x)^{(\frac{n}{2}-\frac{n}{r}+2)}$
 and$\ w_{2}(x)=R(x)^{(3+8k)}$ if we work only with functions
 and $w_{2}(x)=R(x)^{(3+12k)}$ for any $p$-forms, $p\geq 1.$
\end{thm}
\quad As a corollary we get, if we are interested in estimates $L^{r}-L^{s},$\ \par 

\begin{cor}
Let $M$ be a complete non compact riemannian manifold of class
 ${\mathcal{C}}^{2}$ without boundary. Let $\alpha >0,\ r\geq
 2$ and $k$ the smallest integer such that, with $\frac{1}{r_{k}}=\frac{1}{2}-\frac{2k}{n},$
 we have $r_{k}\geq r.$ For any $\omega \in L^{r}(\lbrack 0,T+\alpha
 \rbrack ,L^{r}_{p}(M,w_{1}))\cap L^{r}(\lbrack 0,T+\alpha \rbrack
 ,L^{2}_{p}(M))$ there is a $u\in L^{r}(\lbrack 0,T\rbrack ,L^{s}_{p}(M,w_{2}))$
 such that $\partial _{t}u+\Delta u=\omega $ and:\par 
\vspace{.5em} \ \ \ \ \ \ \ \ \ \ \ \  $\displaystyle {\left\Vert{\partial
 _{t}u}\right\Vert}_{L^{r}(\lbrack 0,T\rbrack ,L^{r}_{p}(M,w_{2}))}+{\left\Vert{u}\right\Vert}_{L^{r}(\lbrack
 0,T\rbrack ,L^{s}_{p}(M,w_{2}))}\leq $\vspace{.5em}  \par 
\vspace{.5em} \ \ \ \ \ \ \ \ \ \ \ \ \ \ \ \ \ \ \ \ \ \ \ \
  $\displaystyle \leq c_{1}{\left\Vert{\omega }\right\Vert}_{L^{r}(\lbrack
 0,T+\alpha \rbrack ,L^{r}_{p}(M,w_{1}))}+c_{2}{\sqrt{T+\alpha
 }}{\left\Vert{\omega }\right\Vert}_{L^{r}(\lbrack 0,T+\alpha
 \rbrack ,L^{2}_{p}(M))}$\vspace{.5em}  \par 
with $\frac{1}{s}=\frac{1}{r}-\frac{2}{n}>0$ and where the weights
 functions are: $w_{1}(x)=R(x)^{(\frac{n}{2}-\frac{n}{r}+2)}$
 and $w_{2}(x)=R(x)^{s(2+(3+8k)/r)}$ if we work only with functions
 and $w_{2}(x)=R(x)^{s(2+(3+12k)/r)}$ for any $p$-forms, $p\geq 1.$
\end{cor}
\quad Proof.\ \par 
By Theorem~\ref{SB32} we have $W_{G}^{m,r}(M,w)\subset W_{G}^{k,s}(M,w'),$
 with $\frac{1}{s}=\frac{1}{r}-\frac{(m-k)}{n}>0.$ So we choose
 $m=2,\ k=0$ and $w=R(x)^{(3+12k)}$ which gives $w'(x)=R(x)^{\nu
 }$ with $\nu =s(2+(3+12k)/r).$ We get then, again with $G=\Lambda
 ^{p}(M)$:\ \par 
\vspace{.5em} \ \ \ \ \ \ \ \ \ \ \ \  $\displaystyle \forall
 u\in W_{p}^{2,r}(M,w),\ {\left\Vert{u}\right\Vert}_{L_{p}^{s}(M,w')}\leq
 C{\left\Vert{u}\right\Vert}_{W_{p}^{2,r}(M,w)}.$\vspace{.5em}  \ \par 
If we work with functions, we have $\nu =s(2+(3+8k)/r).$ \ \par 
Putting this in Theorem~\ref{SB31}, this finishes the proof of
 the corollary. $\hfill\blacksquare $\ \par 

\subsection{Global Gaffney type inequality in $L^{r}.$}
\quad Let $B:=B(x,R)$ be a $(1,\epsilon )$-admissible ball in the complete
 riemannian manifold $(M,g)$ and set $B^{1}:=B(x,R/2).$ Let $\omega
 $ be a $p$-form in $M.$ We have the local $L^{r}$ Gaffney's
 inequality by Corollary~\ref{SC25}:\ \par 
\vspace{.5em} \ \ \ \ \ \ \ \ \ \ \ \  $\displaystyle {\left\Vert{\nabla
 _{M}\omega }\right\Vert}_{L^{r}(B^{1})}\leq C({\left\Vert{d(\omega
 )}\right\Vert}_{L^{r}(B)}+{\left\Vert{d^{*}(\omega )}\right\Vert}_{L^{r}(B)})+cR^{-1}{\left\Vert{\omega
 }\right\Vert}_{L^{r}(B)},$\vspace{.5em}  \ \par 
the constants $c,C$ depending only on $n,\ r,\ p,\ \epsilon .$
 Because $R\leq 1,$ we get, with another constant $c$:\ \par 
\vspace{.5em} \ \ \ \ \ \ \ \ \ \ \ \  $\displaystyle {\left\Vert{\omega
 }\right\Vert}_{W^{1,r}(B^{1})}\leq {\left\Vert{\nabla _{M}\omega
 }\right\Vert}_{L^{r}(B^{1})}+{\left\Vert{\omega }\right\Vert}_{L^{r}(B^{1})}\leq
 $\vspace{.5em}  \ \par 
\vspace{.5em} \ \ \ \ \ \ \ \ \ \ \ \ \ \ \ \ \ \ \ \ \ \ \ \
  $\displaystyle \leq C({\left\Vert{d(\omega )}\right\Vert}_{L^{r}(B)}+{\left\Vert{d^{*}(\omega
 )}\right\Vert}_{L^{r}(B)})+cR^{-1}{\left\Vert{\omega }\right\Vert}_{L^{r}(B)},$\vspace{.5em}
  \ \par 
\quad To globalise this, we proceed as above by summing over the balls
 $B$ and $B^{1}.$ Raising to the power $r$ and adding, we get\ \par 
\vspace{.5em} \ \ \ \ \ \ \ \ \ \ \ \  $\displaystyle \sum_{x\in
 {\mathcal{D}}(\epsilon )}{{\left\Vert{\omega }\right\Vert}^{r}_{W^{1,r}(B^{1})}}\leq
 C\sum_{x\in {\mathcal{D}}(\epsilon )}{({\left\Vert{d(\omega
 )}\right\Vert}^{r}_{L^{r}(B)}+{\left\Vert{d^{*}(\omega )}\right\Vert}^{r}_{L^{r}(B)}+cR^{-r}{\left\Vert{\omega
 }\right\Vert}^{r}_{L^{r}(B)})}$\vspace{.5em}  \ \par 
\quad As above we have\ \par 
\vspace{.5em} \ \ \ \ \ \ \ \ \ \ \ \  $\displaystyle {\left\Vert{\omega
 }\right\Vert}^{r}_{W^{1,r}(M)}\leq \sum_{x\in {\mathcal{D}}(\epsilon
 )}{{\left\Vert{\omega }\right\Vert}^{r}_{W^{1,r}(B^{1})}}$\vspace{.5em}
  \ \par 
and, because the overlap of the covering $\lbrace B(x,R_{\epsilon
 }(x)),\ x\in {\mathcal{D}}(\epsilon )\rbrace $ is bounded by
 $T_{1}$ by Corollary~\ref{SC26}, we get\ \par 
\vspace{.5em} \ \ \ \ \ \ \ \ \ \ \ \  $\displaystyle \sum_{x\in
 {\mathcal{D}}(\epsilon )}{({\left\Vert{d(\omega )}\right\Vert}^{r}_{L^{r}(B)}+{\left\Vert{d^{*}(\omega
 )}\right\Vert}^{r}_{L^{r}(B)})}\leq T_{1}({\left\Vert{d(\omega
 )}\right\Vert}^{r}_{L^{r}(M)}+{\left\Vert{d^{*}(\omega )}\right\Vert}^{r}_{L^{r}(M)})$\vspace{.5em}
  \ \par 
and\ \par 
\vspace{.5em} \ \ \ \ \ \ \ \ \ \ \ \  $\displaystyle \sum_{x\in
 {\mathcal{D}}(\epsilon )}{R_{\epsilon }^{-r}(x){\left\Vert{\omega
 }\right\Vert}^{r}_{L^{r}(B(x,R_{\epsilon }(x)))})}\leq T_{1}{\left\Vert{\omega
 }\right\Vert}^{r}_{L^{r}(M,w)},$\vspace{.5em}  \ \par 
where $w(x)$ is the weight $w(x):=R^{-r}_{\epsilon }(x).$ So
 we proved the global Gaffney's type inequality with weight:\ \par 

\begin{thm}
~\label{SB44}Let $(M,g)$ be a complete riemannian manifold. Let
 $r\geq 1$ and $w(x):=R(x)^{-r}.$ Let $\omega $ be a $p$-form
 in $M.$ We have:\par 
\vspace{.5em} \ \ \ \ \ \ \ \ \ \ \ \  $\displaystyle {\left\Vert{\omega
 }\right\Vert}_{W_{p}^{1,r}(M)}\leq C({\left\Vert{d(\omega )}\right\Vert}_{L_{p+1}^{r}(M)}+{\left\Vert{d^{*}(\omega
 )}\right\Vert}_{L_{p-1}^{r}(M)}+{\left\Vert{\omega }\right\Vert}_{L_{p}^{r}(M,w)}).$\vspace{.5em}
  
\end{thm}
\quad So using Theorem~\ref{SB32} with $\gamma =0,\ m=1,\ k=0,\ \frac{1}{s}=\frac{1}{r}-\frac{1}{n}>0,\
 w'(x)=R_{\epsilon }(x)^{2s},$ plus Theorem~\ref{SB44}, we get\ \par 

\begin{cor}
~\label{SC29}Let $(M,g)$ be a complete riemannian manifold. Let
 $r\geq 1$ and $w(x):=R(x)^{-r}.$ Let $\omega $ be a $p$-form
 in $M.$ We have:\par 
\vspace{.5em} \ \ \ \ \ \ \ \ \ \ \ \  $\displaystyle {\left\Vert{\omega
 }\right\Vert}_{L_{p}^{s}(M,w')}\leq C({\left\Vert{d(\omega )}\right\Vert}_{L_{p+1}^{r}(M)}+{\left\Vert{d^{*}(\omega
 )}\right\Vert}_{L_{p-1}^{r}(M)}+{\left\Vert{\omega }\right\Vert}_{L_{p}^{r}(M,w)})$\vspace{.5em}
  \par 
with $\frac{1}{s}=\frac{1}{r}-\frac{1}{n}>0,\ w'(x)=R_{\epsilon }(x)^{2s}.$
\end{cor}
\quad This corollary is a kind of N. Lohou\'e's result~\cite{Lohoue85}
 with weights and without any geometric conditions on the riemannian
 manifold $(M,g).$\ \par 
\ \par 

\section{Applications.~\label{SB41}}
\quad We shall give some examples where we have classical estimates
 using that $\forall x\in M,\ R_{\epsilon }(x)\geq \delta ,$
 via~\cite[Corollary, p. 7] {HebeyHerzlich97} (see also Theorem
 1.3 in the book by Hebey~\cite{Hebey96}).\ \par 
\quad We have:\ \par 

\begin{cor}
~\label{pL32}Let $(M,g)$ be a complete riemannian manifold. If
 the injectivity radius verifies $r_{inj}(x)\geq i>0$ and the
 Ricci curvature verifies $Rc_{(M,g)}(x)\geq \lambda g_{x}$ for
 some $\lambda \in {\mathbb{R}}$ and all $x\in M,$ then there
 exists a positive constant $\delta >0,$ depending only on $n,\epsilon
 ,\lambda ,\alpha ,i$ such that for any $x\in M,\ r_{H}(1+\epsilon
 ,0,\alpha )(x)\geq \delta .$ \par 
Hence $R_{0,\epsilon }(x)\geq \delta $ all $x\in M.$\par 
\quad If we have $\left\vert{Rc_{(M,g)}(x)}\right\vert \leq C$ for
 all $x\in M,$ then there exists a positive constant $\delta
 >0,$ depending only on $n,\epsilon ,C,\alpha ,i$ and  $c,$ such
 that for any $x\in M,\ r_{H}(1+\epsilon ,1,\alpha )(x)\geq \delta .$\par 
Hence $R_{1,\epsilon }(x)\geq \delta $ all $x\in M.$
\end{cor}
\quad Proof.\ \par 
The Theorem of Hebey and Herzlich gives that, under these hypotheses,
 for any $\alpha \in (0,1)$ that $\forall x\in M,\ r_{H}(1+\epsilon
 ,0,\alpha )(x)\geq \delta .$\ \par 
\quad Now recall that $r_{H}(1+\epsilon ,0,\alpha )(x)$ is the sup
 of $S$ such that, in an harmonic coordinates patch,\ \par 
\quad 1) $(1-\epsilon )\delta _{ij}\leq g_{ij}\leq (1+\epsilon )\delta
 _{ij}$ in $B(x,S)$ as bilinear forms,\ \par 
\quad 2) $S^{\alpha }(x)\sup \ _{i,j=1,...,n,\ y\neq z\in B(x,S)}\left\vert{\frac{g_{ij}(y)-g_{ij}(z)}{d_{g}(y,z)^{\alpha
 }}}\right\vert \leq \epsilon .$\ \par 
So when we take the sup for $R_{0,\epsilon }(x)$ on \emph{any}
 smooth coordinates patch we get $r_{H}(1+\epsilon ,0,\alpha
 )(x)\leq R_{0,\epsilon }(x).$\ \par 
\quad The same way we get $r_{H}(1+\epsilon ,1,\alpha )(x)\leq R_{1,\epsilon
 }(x).$\ \par 
\quad The proof is complete. $\hfill\blacksquare $\ \par 
\ \par 
\quad As a corollary we retrieve Corollary 3.19, p. 38, in~\cite{Hebey96}:\ \par 

\begin{cor}
~\label{SB35}Let $(M,g)$ be a complete riemannian manifold. If
 we have the injectivity radius bounded below and the Ricci curvature
 verifying $Rc(x)\geq \lambda g_{x}$ for some $\lambda \in {\mathbb{R}}$
 and all $x\in M$ then the Sobolev embeddings for functions are
 valid in $(M,g).$
\end{cor}

\begin{rem}
Because the proof of the Theorem of Hebey and Herzlich does not
 use the Theorem of Varopoulos, we get here a different proof
 of Corollary~\ref{SB35}.
\end{rem}
\quad We get also a Sobolev embedding for sections of $G.$\ \par 

\begin{cor}
Let $(M,g)$ be a complete riemannian manifold. Let $G:=(H,\pi
 ,M)$ be a complex ${\mathcal{C}}^{m}$ adapted vector bundle
 over $M.$ If $M$ has $0$-order weak bounded geometry, then the
 "classical" Sobolev embeddings for sections of $G$ are valid in $(M,g).$
\end{cor}
\quad As already said in the introduction, this improves a very well
 known Theorem by M. Cantor~\cite{Cantor74}. Also our proof is
 completely different.\ \par 
\ \par 
\quad We also have a global Gaffney's type inequality:\ \par 

\begin{cor}
Let $(M,g)$ be a complete riemannian manifold with a $0$-order
 weak bounded geometry, then the global Gaffney's type inequality
 in $L^{r}$ is valid:\par 
\vspace{.5em} \ \ \ \ \ \ \ \ \ \ \ \  $\displaystyle {\left\Vert{\omega
 }\right\Vert}_{W_{p}^{1,r}(M)}\leq C({\left\Vert{d(\omega )}\right\Vert}_{L_{p+1}^{r}(M)}+{\left\Vert{d^{*}(\omega
 )}\right\Vert}_{L_{p-1}^{r}(M)}+{\left\Vert{\omega }\right\Vert}_{L_{p}^{r}(M)}).$\vspace{.5em}
  
\end{cor}
\quad Also Corollary~\ref{SC29} take the form:\ \par 

\begin{cor}
Let $(M,g)$ be a complete riemannian manifold with a $0$-order
 weak bounded geometry. Let $r\geq 1$and let $\omega $ be a $p$-form
 in $M.$ We have:\par 
\vspace{.5em} \ \ \ \ \ \ \ \ \ \ \ \  $\displaystyle {\left\Vert{\omega
 }\right\Vert}_{L_{p}^{s}(M)}\leq C({\left\Vert{d(\omega )}\right\Vert}_{L_{p+1}^{r}(M)}+{\left\Vert{d^{*}(\omega
 )}\right\Vert}_{L_{p-1}^{r}(M)}+{\left\Vert{\omega }\right\Vert}_{L_{p}^{r}(M)})$\vspace{.5em}
  \par 
with $\frac{1}{s}=\frac{1}{r}-\frac{1}{n}>0.$
\end{cor}
\quad N. Lohou\'e~\cite{Lohoue85} proved the same result under the
 stronger hypothesis that $(M,g)$ has a $2$-order bounded geometry
 plus some other hypotheses on the laplacian and the range of $r.$\ \par 

\subsection{Weak bounded geometry}
\quad The Corollary~\ref{pL32} gives, in particular, that if $(M,g)$
 has a $0$-order weak bounded geometry, then we get that $\forall
 x\in M,\ R_{(1,\epsilon )}(x)\geq \delta .$ Hence we get the
 improved classical Sobolev embedding Theorems:\ \par 

\begin{thm}
Let $(M,g)$ be a complete riemannian manifold. Let $G:=(H,\pi
 ,M)$ be a complex smooth adapted vector bundle over $M.$ Suppose
 $(M,g)$ has a $0$-order weak bounded geometry. Let $0\leq k<m$
 and $1/s=1/r-(m-k)/n.$ Let $u\in W^{m,r}_{G}(M).$ Then we have
 $u\in W^{k,s}_{G}(M)$ with the control:\par 
\vspace{.5em} \ \ \ \ \ \ \ \ \ \ \ \  $\displaystyle {\left\Vert{u}\right\Vert}^{s}_{W^{k,s}_{G}(M)}\leq
 C{\left\Vert{u}\right\Vert}^{r}_{W^{m,r}_{G}(M)}.$\vspace{.5em}  
\end{thm}
\quad We shall give some examples of such a situation.\ \par 

\subsection{Examples of manifolds of bounded geometry.}
\quad These examples are precisely taken from ~\cite[p. 40]{Eldering12},
 where there is more explanations.\ \par 
\quad $\bullet $  Euclidean space with the standard metric trivially
 has bounded geometry.
\ \par 
\ \par 
\quad $\bullet $ A smooth, compact Riemannian manifold $M$ has bounded
 geometry as well;
 both the injectivity radius and the curvature
 including derivatives are continuous
 functions, so these attain
 their finite minima and maxima, respectively,
on $M.$ If $M\in
 {\mathcal{C}}^{m+2},$ then it has bounded geometry of order $m.$\ \par 
\quad $\bullet $ Non compact, smooth Riemannian manifolds that possess
 a transitive group
 of isomorphisms (such as the hyperbolic
 spaces ${\mathbb{H}}^{n}$) have $m$-order bounded geometry since
 the
 finite injectivity radius and curvature estimates at any
 single point translate to
 a uniform estimate for all points
 under isomorphisms.\ \par 
\ \par 
\quad More manifolds of bounded geometry can be constructed with these
 basic building
 blocks in the following ways.\ \par 
\quad $\bullet $
 The product of a finite number of manifolds of bounded
 geometry again has
 bounded geometry, since the direct sum structure
 of the metric is inherited
 by the exponential map and curvature.\ \par 
\quad $\bullet $ If we take a finite connected sum of manifolds with
 bounded geometry such
 that the gluing modifications are smooth
 and contained in a compact set, then
 the resulting manifold
 has bounded geometry again.\ \par 

\subsection{Compact riemannian manifold with smooth boundary.~\label{SB42}}
\quad We have seen that a compact riemannian manifold of class ${\mathcal{C}}^{2}$
 without boundary has a $0$-order bounded geometry, so the Sobolev
 embeddings are valid on it. We shall deduce that they are also
 valid in the case of a compact riemannian manifold with a ${\mathcal{C}}^{\infty
 }$ smooth boundary. But first we have, using that a compact
 riemannian manifold has $0$-order bounded geometry:\ \par 

\begin{cor}
~\label{SB33}Let $(M,g)$ be a compact riemannian manifold without
 boundary. Let $G:=(H,\pi ,M)$ be an adapted complex smooth vector
 bundle over $M.$\par 
Then the Sobolev embeddings for sections of $G$ are valid in
 $(M,g).$ Precisely we have: $W_{G}^{m,r}(M)$ is embedded in
 $W_{G}^{k,s}(M),$ with $\frac{1}{s}=\frac{1}{r}-\frac{(m-k)}{n}>0$ and:\par 
\vspace{.5em} \ \ \ \ \ \ \ \ \ \ \ \  $\displaystyle \forall
 u\in W_{G}^{m,r}(M),\ {\left\Vert{u}\right\Vert}_{W_{G}^{k,s}(M)}\leq
 C{\left\Vert{u}\right\Vert}_{W_{G}^{m,r}(M)}.$\vspace{.5em}  
\end{cor}
\quad Now let $M$ be a ${\mathcal{C}}^{\infty }$ smooth connected compact
 riemannian manifold with a ${\mathcal{C}}^{\infty }$ smooth
 boundary $\partial M.$ We want to show how the results in case
 of a compact boundary-less manifold apply to this case.\ \par 
\quad A classical way to get rid of an "annoying boundary" of a manifold
 is to use its "double". For instance: Duff~\cite{Duff52}, H\"ormander~\cite[p.
 257]{Hormand94}. \ \par 
Here we copy the following construction from~\cite[Appendix B]{GuneysuPigola}.\
 \par 
\quad The "Riemannian double" $\Gamma :=\Gamma (M)$ of $M,$ obtained
 by gluing two copies, $M$ and $M_{2},$ of $M$ along $\partial
 M,$ is a compact Riemannian manifold without boundary. Moreover,
 by its very construction, it is always possible to assume that
 $\displaystyle \Gamma $ contains an isometric copy of the original
 manifold $M.$ We shall also write $M$ for this isometric copy
 to ease notation.\ \par 
\ \par 
\quad We take $u\in W_{G}^{m,r}(M)$ and we want to show that $u\in
 W_{G}^{k,s}(M),$ with $\frac{1}{s}=\frac{1}{r}-\frac{(m-k)}{n}>0.$\ \par 
\quad We shall suppose that $G$ extends smoothly to $\Gamma ,$ i.e.
 the connexion is smooth and still is a metric connexion on $\Gamma
 ,$ and the scalar product also is smooth in $\Gamma .$ For instance
 this is the case if $G=\Lambda ^{p}(M),$ the bundle of $p$-forms in $M.$\ \par 
\quad The Seeley Theorem~\cite{Seeley64} in the version of Lions~\cite{Lions64},
 tells us that any function $f\in W^{m,r}(M)$ can be extended
 to $\Gamma $ as $\tilde f\in W^{m,r}(\Gamma ).$ By use of a
 finite covering of $M$ by balls $B(x,R(x))$ with center $x\in
 M$ and trivializing the bundle $G,$ and an associated partition
 of unity, this result of Seeley can be made valid to a section
 $u\in W_{G}^{m,r}(M).$ So we have an extension $\tilde u\in
 W_{G}^{m,r}(\Gamma ).$\ \par 
\quad Using Corollary~\ref{SB33} we get that $\tilde u\in W_{G}^{k,s}(\Gamma
 ),$ with $\frac{1}{s}=\frac{1}{r}-\frac{(m-k)}{n}>0$ and:\ \par 
\vspace{.5em} \ \ \ \ \ \ \ \ \ \ \ \  $\displaystyle \forall
 \tilde u\in W_{G}^{m,r}(\Gamma ),\ {\left\Vert{\tilde u}\right\Vert}_{W_{G}^{k,s}(\Gamma
 )}\leq C{\left\Vert{\tilde u}\right\Vert}_{W_{G}^{m,r}(\Gamma
 )}.$\vspace{.5em}  \ \par 
Hence, restricting $\tilde u$ to $M$ we get, a fortiori, $u\in
 W_{G}^{k,s}(M),$ and:\ \par 
\vspace{.5em} \ \ \ \ \ \ \ \ \ \ \ \  $\displaystyle \forall
 u\in W_{G}^{m,r}(M),\ {\left\Vert{u}\right\Vert}_{W_{G}^{k,s}(M)}\leq
 C{\left\Vert{u}\right\Vert}_{W_{G}^{m,r}(M)}.$\vspace{.5em}  \ \par 
\quad So we proved:\ \par 

\begin{thm}
Let $(M,g)$ be a compact riemannian manifold with a smooth boundary.
 Let $G:=(H,\pi ,M)$ be a complex vector bundle over $M,$ which
 admit a smooth adapted extension to a "double" manifold $\Gamma
 .$ Then the Sobolev embeddings for sections of $G$ are valid
 in $(M,g).$ Precisely we have: $W_{G}^{m,r}(M)$ is embedded
 in $W_{G}^{k,s}(M),$ with $\frac{1}{s}=\frac{1}{r}-\frac{(m-k)}{n}>0$ and:\par 
\vspace{.5em} \ \ \ \ \ \ \ \ \ \ \ \  $\displaystyle \forall
 u\in W_{G}^{m,r}(M),\ {\left\Vert{u}\right\Vert}_{W_{G}^{k,s}(M)}\leq
 C{\left\Vert{u}\right\Vert}_{W_{G}^{m,r}(M)}.$\vspace{.5em}  
\end{thm}

\subsection{Hyperbolic manifolds.~\label{SB43}}
\quad These are manifolds such that the sectional curvature $K_{M}$
 is constantly $-1.$ For them we have first that the Ricci curvature
 is bounded.\ \par 

\begin{lem}
~\label{BG0} Let $(M,g)$ be a complete Riemannian manifold such
 that $H\leq K_{M}\leq K$ for constants
 $H,K\in {\mathbb{R}}.$
 Then we have that ${\left\Vert{Rc}\right\Vert}_{\infty }\leq
 \max (\left\vert{H}\right\vert ,\left\vert{K}\right\vert ).$
\end{lem}
\quad Proof.\ \par 
Take a tangent vector $v\in T_{x}M,$ with $\left\vert{v}\right\vert
 =1,$ we obtain the Ricci curvature $Rc$ of $v$ at $x$ by extending
 $v=v_{n}$ to an orthonormal basis $v_{1},...,v_{n}.$ Then we
 compute the Ricci curvature along $v$:\ \par 
\vspace{.5em} \ \ \ \ \ \ \ \ \ \ \ \  $\displaystyle Rc_{x}(v)=\frac{1}{n-1}\sum_{j=1}^{n-1}{{\left\langle{R_{x}(v,v_{j})v,v_{j}}\right\rangle}}.$\vspace{.5em}
  \ \par 
where $R$ denotes the Riemannian curvature tensor. On the other
 hand, the sectional curvature $K_{M,x}(v,v_{j})$ for $j<n$ is
 given by (remember the $v_{j}$ are orthonormal):\ \par 
\vspace{.5em} \ \ \ \ \ \ \ \ \ \ \ \  $\displaystyle K_{M,x}(v,v_{j})={\left\langle{R_{x}(v,v_{j})v,v_{j}}\right\rangle}$\vspace{.5em}
  \ \par 
So we get\ \par 
\vspace{.5em} \ \ \ \ \ \ \ \ \ \ \ \  $\displaystyle Rc_{x}(v)=\frac{1}{n-1}\sum_{j=1}^{n-1}{{\left\langle{R_{x}(v,v_{j})v,v_{j}}\right\rangle}}=\frac{1}{n-1}\sum_{j=1}^{n-1}{K_{M,x}(v,v_{j})}$\vspace{.5em}
  \ \par 
hence $H\leq Rc_{x}(v)\leq K\Rightarrow \left\vert{Rc_{x}}\right\vert
 \leq \max (\left\vert{H}\right\vert ,\left\vert{K}\right\vert ).$\ \par 
\quad To get estimates on the Ricci tensor $Rc_{(M,g)}(x)(u,v),$ we
 notice that $Rc_{x}(v)=Rc_{(M,g)}(x)(v,v)$ and we get the estimates
 by polarisation. $\hfill\blacksquare $\ \par 
\ \par 
\quad To get that the injectivity radius $r_{inj}(x)$ is bounded below
 we shall use a Theorem by Cheeger, Gromov
 and Taylor~\cite{CheegerGroTay82}:\
 \par 

\begin{thm}
~\label{pL33}Let $(M,g)$ be a complete Riemannian manifold such
 that $K_{M}\leq K$ for constants
 $K\in {\mathbb{R}}.$ Let $0<r<\frac{\pi
 }{4{\sqrt{K}}}$ if $K>0$ and $r\in (0,\infty )$ if $K\leq 0.$
 Then the injectivity radius $r_{inj}(x)$ at $x$ satisfies\par 
\vspace{.5em} \ \ \ \ \ \ \ \ \ \ \ \  $\displaystyle r_{inj}(x)\geq
 r\frac{\mathrm{V}\mathrm{o}\mathrm{l}(B_{M}(x,r))}{\mathrm{V}\mathrm{o}\mathrm{l}(B_{M}(x,r))+\mathrm{V}\mathrm{o}\mathrm{l}(B_{T_{x}M}(0,2r))},$\vspace{.5em}
  \par 
where $B_{T_{x}M}(0,2r))$ denotes the volume of the ball of radius
 $2r$ in $T_{x}M,$ where both the
 volume and the distance function
 are defined using the metric $g^{*}:=\exp _{p}^{*}g$ i.e. the
 pull-back of
 the metric $g$ to $T_{x}M$ via the exponential map.
\end{thm}
\quad This Theorem leads to the definition:\ \par 

\begin{defin}
Let $(M,g)$ be a Riemannian manifold. We shall say that it has
 the {\bf lifted doubling property} if we have:\par 
\vspace{.5em} \ \ \ \ \ \ \ \ \ \ \ \  $\displaystyle (LDP)\
 \ \ \ \ \ \ \ \ \exists \beta ,\gamma >0::\forall x\in M,\ \exists
 r\geq \beta ,\ \mathrm{V}\mathrm{o}\mathrm{l}(B_{T_{x}M}(0,2r))\leq
 \gamma \mathrm{V}\mathrm{o}\mathrm{l}(B_{M}(x,r)),$\vspace{.5em}  \par 
where $B_{T_{x}M}(0,2r))$ denotes the volume of the ball of radius
 $2r$ in $T_{x}M,$ and both the
 volume and the distance function
 are defined on $T_{x}M$ using the metric $g^{*}:=\exp _{p}^{*}g$
 i.e. the pull-back of
 the metric $g$ to $T_{x}M$ via the exponential map.
\end{defin}
Hence we get:\ \par 

\begin{cor}
~\label{pL44}Let $(M,g)$ be a complete Riemannian manifold such
 that $K_{M}\leq K$ for a constant
 $K\in {\mathbb{R}}.$ For
 instance an hyperbolic manifold.  Suppose moreover that $(M,g),$
 has the lifted doubling property. Then:\par 
\vspace{.5em} \ \ \ \ \ \ \ \ \ \ \ \  $\displaystyle \forall
 x\in M,\ r_{inj}(x)\geq \frac{\beta }{1+\gamma }.$\vspace{.5em}  
\end{cor}
\quad Proof.\ \par 
By the (LDP) we get, for a $r\geq \beta ,$\ \par 
\vspace{.5em} \ \ \ \ \ \ \ \ \ \ \ \  $\displaystyle \mathrm{V}\mathrm{o}\mathrm{l}(B_{T_{x}M}(0,2r))\leq
 \gamma \mathrm{V}\mathrm{o}\mathrm{l}(B_{M}(x,r)).$\vspace{.5em}  \ \par 
We apply Theorem~\ref{pL33} of Cheeger, Gromov
 and Taylor to get\ \par 
\vspace{.5em} \ \ \ \ \ \ \ \ \ \ \ \  $\displaystyle r_{inj}(x)\geq
 r\frac{\mathrm{V}\mathrm{o}\mathrm{l}(B_{M}(x,r))}{\mathrm{V}\mathrm{o}\mathrm{l}(B_{M}(x,r))+\mathrm{V}\mathrm{o}\mathrm{l}(B_{T_{x}M}(0,2r))}.$\vspace{.5em}
  \ \par 
So\ \par 
\vspace{.5em} \ \ \ \ \ \ \ \ \ \ \ \  $\displaystyle \frac{\mathrm{V}\mathrm{o}\mathrm{l}(B_{M}(x,r))}{\mathrm{V}\mathrm{o}\mathrm{l}(B_{M}(x,r))+\mathrm{V}\mathrm{o}\mathrm{l}(B_{T_{x}M}(0,2r))}\geq
 \frac{1}{1+\gamma }$\vspace{.5em}  \ \par 
hence, because $r\geq \beta ,$ we get the result. $\hfill\blacksquare $\ \par 
\ \par 
\quad So finally we get\ \par 

\begin{cor}
Let $(M,g)$ be a complete Riemannian manifold such that $H\leq
 K_{M}\leq K$ for constants
 $H,K\in {\mathbb{R}},$ where $K_{M}$
 is the sectional curvature of $M.$ Let $G:=(H,\pi ,M)$ be an
 adapted complex smooth vector bundle over $M.$ Suppose moreover
 that  $(M,g)$ has the lifted doubling property. Then $\exists
 \delta >0,\ \forall x\in M,\ R_{1,\epsilon }(x)\geq \delta .$\par 
\quad This implies that the Sobolev embeddings are valid for sections
 of $G$ in that case.
\end{cor}
\quad Proof.\ \par 
Because the Lemma~\ref{BG0} gives that the Ricci curvature is
 bounded. The  Corollary~\ref{pL44} gives $\forall x\in M,\ r_{inj}(x)\geq
 \frac{\beta }{1+\gamma }.$ So the Corollary~\ref{pL32} completes
 the proof. $\hfill\blacksquare $\ \par 

\begin{rem}
In the case the hyperbolic manifold $(M,g)$ is simply connected,
 then by the Hadamard Theorem~\cite[Theorem 3.1, p. 149]{DoCarmo93},
 we get that the injectivity radius is $\infty ,$ so we have
 also the classical embedding Theorems in this case, even for sections of $G.$
\end{rem}
\ \par 

\bibliographystyle{/usr/local/texlive/2017/texmf-dist/bibtex/bst/base/apalike}

\begin{thebibliography}{}

\bibitem[Amar, 2018a]{ellipticEq18}
Amar, E. (2018a).
\newblock The {LIR} method. ${L}^r$ solutions of elliptic equation in a
  complete riemannian manifold.
\newblock {\em J. Geometric Analysis}.
\newblock To appear. DOI : 10.1007/s12220-018-0086-3.

\bibitem[Amar, 2018b]{SobPar18}
Amar, E. (2018b).
\newblock Sobolev solutions of parabolic equation in a complete riemannian
  manifold.
\newblock {\em Arxiv}.
\newblock arXiv:1812.O4411.

\bibitem[Aubin, 1982]{Aubin82}
Aubin, T. (1982).
\newblock {\em Nonlinear Analysis on Manifolds. Monge-Amp\'ere Equations},
  volume 252.
\newblock Springer-Verlag, New York.

\bibitem[Cantor, 1974]{Cantor74}
Cantor, M. (1974).
\newblock Sobolev inequalities for riemannian bundles.
\newblock {\em Bull Am. Math. Soc.}, 80:239--243.

\bibitem[Carron, 1994]{Carron94}
Carron, G. (1994).
\newblock In\'egalit\'es isop\'erim\'etriques sur les vari\'et\'es
  riemanniennes.
\newblock Master's thesis, Universit\'e Joseph Fourier, Grenoble.
\newblock These de Doctorat.

\bibitem[Cheeger et~al., 1982]{CheegerGroTay82}
Cheeger, J., Gromov, M., and Taylor, M. (1982).
\newblock Finite propagation speed, kernel estimates for functions of the
  {L}aplace operator, and the geometry of complete riemannian manifolds,.
\newblock {\em J. Differential Geom.}, 17:15--53.

\bibitem[Coulhon and Saloff-Coste, 1993]{CoulSal93}
Coulhon, T. and Saloff-Coste, L. (1993).
\newblock Isop\'erim\'etrie pour les groupes et les vari\'et\'es.
\newblock {\em Revista Mathematica Iberoamericana}, 9:293--314.

\bibitem[do~Carmo, 1993]{DoCarmo93}
do~Carmo, M.~P. (1993).
\newblock {\em Riemannian geometry}.
\newblock Mathematics. Birkh\"auser Boston.

\bibitem[Duff, 1952]{Duff52}
Duff, G. (1952).
\newblock Differential forms in manifolds with boundary.
\newblock {\em Ann. of Math.}, 56:115--127.

\bibitem[Eldering, 2012]{Eldering12}
Eldering, J. (2012).
\newblock Persistence of noncompact normally hyperbolic invariant manifolds in
  bounded geometry.
\newblock {\em Ph. d. Thesis}.
\newblock Ph. d. Thesis Universiteit Utrecht.

\bibitem[Evans and Gariepy, 1992]{EvGar92}
Evans, L.~C. and Gariepy, R.~F. (1992).
\newblock {\em Measure theory and fine properties of functions.}
\newblock Studies in Advanced Mathematics. CRC Press, Boca Raton.

\bibitem[Guneysu and Pigola, 2015]{GuneysuPigola}
Guneysu, B. and Pigola, S. (2015).
\newblock Calderon-{Z}ygmund inequality and {S}obolev spaces on noncompact
  riemannian manifolds.
\newblock {\em Advances in Mathematics}, 281:353--393.

\bibitem[Hebey, 1996]{Hebey96}
Hebey, E. (1996).
\newblock {\em Sobolev spaces on Riemannian manifolds.}, volume 1635 of {\em
  Lecture Notes in Mathematics}.
\newblock Springer-Verlag, Berlin.

\bibitem[Hebey and Herzlich, 1997]{HebeyHerzlich97}
Hebey, E. and Herzlich, M. (1997).
\newblock Harmonic coordinates, harmonic radius and convergence of riemannian
  manifolds.
\newblock {\em Rend. Mat. Appl. (7) 17 (1997), no. 4, 569-605 (1998)},
  17(4):569--605.

\bibitem[H\"ormander, 1994]{Hormand94}
H\"ormander, L. (1994).
\newblock {\em {The Analysis of Linear Partial Differential Operators III}},
  volume 274 of {\em Grundlehren der mathematischen Wissenschften}.
\newblock Springer.

\bibitem[Lions, 1964]{Lions64}
Lions, J.~L. (1964).
\newblock Extension of ${C}^\infty$ functions defined in a half space.
\newblock {\em Math. Reviews, MR0165392 (29 \#2676)}.

\bibitem[Lohou\'e, 1985]{Lohoue85}
Lohou\'e, N. (1985).
\newblock In\'egalit\'es de {S}obolev pour les formes diff\'erentielles sur une
  vari\'et\'e riemannienne.
\newblock {\em C. R. Acad. Sci. Paris}, 301(6):277--280.

\bibitem[Mazzucato and Nistor, 2006]{Nistor06}
Mazzucato, A.-L. and Nistor, V. (2006).
\newblock Mapping properties of heat kernels, maximal regularity, and
  semi-linear parabolic equations on noncompact manifolds.
\newblock {\em Journal of Hyperbolic Differential Equations}, 3(4):599--629.

\bibitem[Scott, 1995]{Scott95}
Scott, C. (1995).
\newblock ${L}^p$ theory of differential forms on manifolds.
\newblock {\em Transactions of the Americain Mathematical Society},
  347(6):2075--2096.

\bibitem[Seeley, 1964]{Seeley64}
Seeley, R. (1964).
\newblock Extension of ${C}^\infty$ functions defined in a half space.
\newblock {\em Proc. Amer. Math. Soc.}, 15:625--626.

\bibitem[Taylor, 2000]{TaylorGD}
Taylor, M.~E. (2000).
\newblock {\em Differential Geometry}.
\newblock Course of the University of North Carolina. University of North
  Carolina.
\newblock www.unc.edu/math/ Faculty/met/diffg.html.

\bibitem[Varopoulos, 1989]{Varopoulos89}
Varopoulos, N. (1989).
\newblock Small time gaussian estimates of heat diffusion kernels.
\newblock {\em Bulletin des Sciences Math\'ematiques.}, 113:253--277.

\bibitem[Voisin, 2002]{Voisin02}
Voisin, C. (2002).
\newblock {\em Th\'eorie de Hodge et g\'eom\'etrie alg\'ebrique complexe.},
  volume~10 of {\em Cours sp\'ecialis\'e}.
\newblock S.M.F.

\end{thebibliography}

\end{document}